\documentclass[leqno]{amsart}
\usepackage{amsthm}
\usepackage{amsfonts}
\usepackage{amsmath, amscd,amssymb,amsbsy,epsf}
    \usepackage{enumerate}
    \usepackage{paralist}
     \usepackage[inline]{enumitem}
\usepackage{bbm, dsfont}
\usepackage{graphicx}
\usepackage{mathrsfs}
\usepackage{verbatim}
\usepackage{inputenc}
\DeclareMathOperator*{\lcm}{lcm}

\def\defining{\overset{\mathbf{def} } = }
\def\R{\boldsymbol{\mathbbm{R} } }
\def\N{\boldsymbol{\mathbbm{N} } }
\def\b{\mathbf{b} }
\def\m{\mathbf{m} }
\def\a{\mathbf{v} }
\def\c{\mathbf{u} }

\def\triang{\mathcal{T}}
\def\triangext{\mathcal{T}_{\mathrm{ext}}}
\def\triangint{\mathcal{T}_{\mathrm{int}}}
\def\triangcon{\mathcal{T}_{\mathrm{nct}}}
\def\vertex{\mathcal{P}}
\def\vertexint{\mathcal{P}_{\mathrm{int}}}
\def\vertexext{\mathcal{P}_{\mathrm{ext}}}

\def\edges{\mathcal{E}_{\scriptscriptstyle \triang}}
\def\edgesint{\mathcal{E}_{\mathrm{int}}}
\def\edgesext{\mathcal{E}_{\mathrm{ext}}}

\def\faces{\mathcal{F}_{\scriptscriptstyle \triang}}
\def\facesint{\mathcal{F}_{\mathrm{int}}}
\def\facesext{\mathcal{F}_{\mathrm{ext}}}

\def\pole{\xi_{\,p}}
\def\Pole{\sigma_{p}}
\def\graph{\mathcal{G}_{\triang}}
\def\edges{\mathcal{E}_{\triang}}

\def\hdelta{h_{\scriptscriptstyle \Delta}}
\def\regdelta{\zeta_{\scriptscriptstyle \Delta}}
\def\adelta{\vartheta_{\scriptscriptstyle \Delta}}
\def\mdelta{\eta_{\scriptscriptstyle \Delta}}
\def\rhodelta{\rho_{\scriptscriptstyle \Delta}}

\def\E{\mathcal{E} }
\def\F{\mathcal{F} }
\def\B{\mathcal{B}}
\def\graph{\mathcal{G}_{\triang}}

\def\rel{\boldsymbol{\sim}}

\def\bspace{\mathcal{R}}
\def\refin{\mathcal{Q}}

\begin{document}

\newtheorem{theorem}{Theorem}[section]
\newtheorem{remark}{Remark}[section]
\newtheorem{corollary}[theorem]{Corollary}
\newtheorem{definition}[theorem]{Definition}
\newtheorem{proposition}[theorem]{Proposition}
\newtheorem{lemma}[theorem]{Lemma}

\oddsidemargin 16.5mm
\evensidemargin 16.5mm

\thispagestyle{plain}

%
%
%
%
\begin{center}

{\large\bf  ON THE GENERATION OF BIPARTITE GRIDS WITH\\[3pt] CONTROLLED REGULARITY FOR 2-D AND 3-D \\SIMPLY CONNECTED DOMAINS
\rule{0mm}{6mm}\renewcommand{\thefootnote}{}
\footnotetext{

\rule{2.4mm}{0mm}Mesh Generation, Shape Measure, Numerical Applications.
}}

\vspace{1cc}
{\large\it Fernando A Morales $\&$ Mauricio A Osorio}

\vspace{1cc}
\parbox{24cc}{{\small

We present a procedure to generate bipartite grids for simply connected domains in 2-D and 3-D of prescribed size and controlled regularity elements. The mesh elements $K$ of the triangulation satisfy $\zeta_{K} \leq C$ where $\zeta_{K}$ is the regularity and $C$ is a constant depending on the shape parameters of the initial mesh. Bipartite grids permit a well-posed mixed-mixed variational formulation of problems such as the porous media flow equation and other linear physical phenomena as well as Galerkin-type discontinuous relaxations.
%

}}
\end{center}

\vspace{1cc}


\vspace{1.5cc}
\section{Introduction}\label{Sec Introduction}
The generation of quality shaped grids for geometric domains is a vast and active research field of mathematics. Defining a grid is one of the key steps in solving partial differential equations with finite element methods. However, the approximation estimates depend not only on the size of the mesh but also on the quality shape of its elements; in addition, other types of numerical difficulties can be introduced in the method due to the quality shape of the elements. It is rather frequent to use triangular and rectangular elements for 2-D domains as well as tetrahedral and hexahedral elements for 3-D domains; the choice is done according to the structure of the modeled problem.  Some other questions in the field are the automatic generation of grids \cite{LentiniPaluszny1}, the generation of polynomial patches for the approximation of two dimensional manifolds, see \cite{Paluszny1}, the generation of structured and/or unstructured grids, the numerical costs and efficiency for grid generation in comparison with the cost of the numerical solutions of the physical problem posed on the grid, the Delaunay tessellation approach for grid generation, used for example for the meshless finite element methods , etc. see \cite{Liseikin,IdelsohnOnate}. On the other hand, the meshfree methods present a totally different approach and try to circumvent the intrinsic difficulties of mesh generation, see \cite{OsorioFrench1}, although other difficulties always arise, see \cite{IdelsohnOnate}.  However, none of these achievements is aimed to fulfill the needs of a particular variational formulation.

In this work, we address the problem of generating bipartite grids and its main motivation is to permit the mixed-mixed variational formulation (presented in \cite{MoralesShow2}) of the partial differential equation. The mixed-mixed formulation is remarkable when modeling problems with interface due to the degrees of freedom introduced in the underlying spaces of functions. Unlike any other formulation for interface problems, there are no linear coupling constraints between function spaces. This feature makes it very attractive for setting discontinuous Galerkin schemes and analyzing problems with micro structure of fractal type. Although the formulation has been developed only for the problem of saturated flow in porous media, it can be easily extended to problems of linear elasticity, heat diffusion, etc. In the problem presented in \cite{MoralesShow2} the domain of analysis $\Omega$ is subdivided in two, namely $\Omega^{1}$ and $\Omega^{2}$ the first uses the pairing $[\mathbf{H}_{\mathrm{div} }, L^{2}]$ and the second uses the pair $[\mathbf{L}^{\!2}, H^{1}]$ for the velocity and pressure respectively. However in that particular case it is clear how to subdivide the domain and set the pairs of modeling spaces of functions. If the formulation is to be used with another subdivision of the domain, namely the one provided by a mesh on $\Omega$, we need to assure that each element with pairing $[\mathbf{H}_{\mathrm{div} }, L^{2}]$ shares boundary of non-negligible measure, only with elements whose pairing is of the type $[\mathbf{L}^{\!2}, H^{1}]$. Therefore, the graph of the mesh has to be bipartite, where two elements of the grid are connected if they have a common interface of non-negligible measure.

Due to the limitations of the technique we restrict our attention to simply connected domains in $\R^{\!2}$ and $\R^{\!3}$. For simplicity it will be also assumed that the domains are polygonal or polyhedral. This is not conceptually far from a more realistic case since it is straightforward to generate a mesh of curved triangles or tetrahedra for a simply connected domain using the gridding of a polygon or polyhedron inscribed in the original domain. Under the hypothesis above the necessary and sufficient condition for triangular and tetrahedral grids to be bipartite will be presented. Additionally, a method for generating grids of arbitrary small size, with elements of bounded regularity will also be introduced. Both, the characterization of bipartite grids and the method to generate them fine and regular, are quite simple and easy to implement. However, it will be clear that the generation method is not optimal; our main goal is to provide the theoretical setting to assure the well-posed mixed-mixed formulation of the problem as well as the convergence of approximate solutions.

Next, we introduce the notation. In the following $\Omega$ denotes an open simply connected polygonal or polyhedral bounded region of $\R^{\!2}$ or $\R^{\!3}$ respectively. Vectors in $\R^{\!2}$ or $\R^{\!3}$ are denoted with bold letters. We write $\#A$ for the cardinal of the set $A$ and $\vert A \vert$ for the length of a segment;  $\vert  A\vert_{i}$ with $i = 2, 3$ stands for the area and volume of the set, depending on the context. Triangles and tetrahedra will typically be denoted with the letters $K, L, M$ and $\Delta$ and since the notation is consistent they shall be seen as elements of the grid or as vertices of its associated graph depending on the context. We write $\triang$ for triangulations or tetrahedral grids of the domain and the characters $\B, \refin, \mathcal{S}$ for refinement processes of the mesh. In the reminder of this section we indicate the minimum background necessary from graph theory \cite{GrossYellen}.

%
%
\subsection{Preliminaries from Graph Theory}
%
%
%
%
\begin{definition} \label{Def walk, path and cycle, length}
Let $G = (V, E)$ be a graph
\\[5pt]
\begin{enumerate*}[label = (\roman*), itemjoin={\\[5pt]\noindent}, mode=unboxed]
\item
A \textbf{walk} in $G$ from vertex $v_{0}$ to vertex $v_{j}$ is an alternating sequence
\begin{equation*}
W \defining \langle v_{0}, e_{1}, v_{1}, e_{2}, \ldots, v_{n-1}, e_{n}, v_{n}\rangle
\end{equation*}
of vertices and arcs such that the endpoints of the edge $e_{i}$ are $v_{i-1}$ and $v_{i}$ for all $i = 1, 2, \ldots, n$.

\item
A \textbf{path} is a walk with no repeated arcs and no repeated vertices.

\item
A walk or path is \textbf{trivial} if it has only one vertex and no arcs.

\item
A \textbf{cycle} is a non-trivial closed path \emph{i.e.} it starts and ends on the same vertex.

\item
The \textbf{length} of a walk, path or cycle is the number of arc-steps in the sequence. We denote it by $\# \langle v_{0}, e_{1}, v_{1}, e_{2}, \ldots, v_{n-1}, e_{n}, v_{n}\rangle$.
\end{enumerate*}
\end{definition}
\begin{definition} A \textbf{cycle graph} $C = (V_{C}, E_{C})$ is a single vertex with a self-loop or a \textbf{simple graph} with $\# V_{C} = \# E_{C}$ that can be drawn so that all its vertices and edges lie on a single circle. A $j$-vertex cycle graph is denoted $C_{j}$.
\end{definition}
\begin{definition}\label{Def bipartite graph}
A \textbf{bipartite graph} is a graph whose vertex set $V$ can be partitioned in two subsets $U, W$ such that each edge of $G$ has one endpoint in $U$ and one endpoint in $W$. The pair $U, W$ is called a (vertex) bipartition of $G$ and $U$ and $W$ are called the bipartition subsets.
\end{definition}
Next we recall a well-known characterization result for bipartite graphs \cite{GrossYellen}
\begin{theorem}\label{Th bipartite graphs characterization}
A graph $G$ is bipartite if and only if it has no cycles of odd length.
\end{theorem}
\begin{remark}\label{Rem language clarification}
   Throughout the present work, we will call vertices the elements of the studied graph as well as the geometric vertices of the triangles or tetrahedra. This must be understood depending on the context.
\end{remark}
%
%
\section{The Two Dimensional Case}
%
%
%
%
We start this section defining the graph associated to the triangulation of a given polygonal domain in $\R^{\! 2}$.
\begin{definition}\label{Def triangulation}
Let $\mathcal{O}\subseteq \R^{\! 2}$ be an open bounded polygonal domain and $\{K: K\in \triang\}$ be a triangulation of $\mathcal{O}$. We denote $\vertex$ 
the set of vertices 
of the triangulation.
%
\\[5pt]
\begin{enumerate*}[label = (\roman*), itemjoin={\\[5pt]\noindent}, mode=unboxed]

\item
The associated graph $\graph(\triang, \edges)$ is defined in the following way. The set of vertices $\triang$ is defined by the set of triangles and there is an arc in $\edges$ between two elements $K, L\in \triang$ if $\vert \partial K\cap \partial L\vert  > 0$ \emph{i.e.} if they share a common side.

\item
An element of the triangulation $\{K:K\in \triang\}$ is said to be isolated if it shares no common side with any other triangle. Equivalently if its degree as vertex of $\graph$ is zero. We denote $\triangcon$ the connected or no isolated, elements of the triangulation.

\item
We say a triangle is INTERIOR if $\vert \partial K\cap \partial \mathcal{O}\vert = 0$ and denote $\triangint$ the set of interior triangles. We say a triangle $K\in \triang$ is EXTERIOR if $\vert \partial K\cap \partial \mathcal{O}\vert>0$ and denote $\triangext$ the set of exterior triangles.

\item
We say a vertex $\xi\in \vertex$ is EXTERIOR if it lies on the boundary of the domain $\partial \mathcal{O}$ and denote the set of exterior vertices by $\vertexext$. We say a vertex $\xi\in \vertex$ of a triangle is INTERIOR if it belongs to the interior of the domain and denote the set of interior vertices by $\vertexint$.

\end{enumerate*}
\end{definition}
\begin{remark}\label{Rem Graph Connectedness}
   Notice that if a polygonal domain $\Omega$ is simply connected its graph $\graph$ has to be connected \emph{i.e.} $\triang = \triangcon$.
\end{remark}
\begin{definition}\label{Def Bipartite triangulation}
Let $\mathcal{O}\subseteq \R^{\! 2}$ be an open bounded polygonal domain and $\{K: K\in \triang\}$ be a triangulation of $\mathcal{O}$. We say a triangulation $\{K : K\in \triang\}$ is BIPARTITE if it can be colored with only two colors. Equivalently if its associated graph $\graph(\triang, \edges)$ is bipartite.
\end{definition}
\begin{definition}\label{Def control points and paths}
Let $\mathcal{O}\subseteq \R^{\!2}$ be an open bounded polygonal domain and $\{K: K\in \triang\}$ be a triangulation of $\mathcal{O}$.
\\[5pt]
\begin{enumerate*}[label = (\roman*), itemjoin={\\[5pt]\noindent}, mode=unboxed]
\item
For each $K\in \triang$ we denote $\b_{K}$ its barycenter or center of gravity.

\item
Given a cycle $C = \langle K_{1}, K_{2}, \ldots, K_{j}, K_{1}\rangle$ in the graph $\graph$ we denote $\gamma_{\scriptscriptstyle C}$ the rectifiable path generated by the sequence of segments $[\b_{\scriptscriptstyle K_{1}}, \b_{\scriptscriptstyle K_{2}}], \ldots [\b_{\scriptscriptstyle K_{j-1}},\b_{\scriptscriptstyle K_{j}}], [\b_{\scriptscriptstyle K_{j}},\b_{\scriptscriptstyle K_{j}}]$ \emph{i.e.} $\gamma_{\scriptscriptstyle C}$ is contained in $\mathcal{O}$.

\item
Given a rectifiable path $\gamma$ and a point $\b$ in $\R^{\!2} - \{\gamma\}$ we denote $n(\gamma, \b)$ the winding number \cite{ConwayComplex}.
\end{enumerate*}
\end{definition}
\begin{remark}\label{Rem Path Cycle}
Let $C$ be the a cycle in the associated graph $\graph$ of a triangulation $\triang$ notice the following
\\[5pt]
\begin{enumerate*}[label = (\roman*), itemjoin={\\[5pt]\noindent}, mode=unboxed]
\item
The path $\gamma_{\scriptscriptstyle C}$ divides the plane in only two connected components. The application $\b\mapsto n(\gamma_{\scriptscriptstyle C}, \b)$ defined on $\R^{\!2} - \{\gamma_{\scriptscriptstyle C}\}$ takes only three values: $0$ on the unbounded connected component and, on the bounded component, $1$ if $\gamma_{\scriptscriptstyle C}$ is counterclockwise oriented or $-1$ if it is clockwise oriented.

\item
Let $K$ be an element of the triangulation and $\{\xi_{i}:1\leq i\leq 3\}$ its three vertices. Due to the definition of the path $\gamma_{\scriptscriptstyle C}$ we know $\{\xi_{i}:1\leq i\leq 3\}\cap \{\gamma_{\scriptscriptstyle C}\} = \emptyset$, hence, the winding number $n(\gamma_{\scriptscriptstyle C}, \xi)$ is well-defined for all $\xi\in \vertex$.
\end{enumerate*}
\end{remark}
%
%
%
%
%
%
\subsection{The Characterization}
%
%
First we focus on a very particular type of triangulation.
\begin{definition}\label{Def radial triangulation}
We say a triangulation $\{K: K\in \triang\}$ of a polygonal domain $\mathcal{O}$ is RADIAL if all its elements are connected and it has only one interior vertex denoted $\pole$ such that $\{K\in \triang: \pole\in \partial K\} = \triang$. In the following we will refer to this interior vertex as POLE.
\end{definition}
\begin{proposition}\label{Th radial mesh}
Let $\{K: K\in \triang\}$ be a radial triangulation of the domain $\Omega$ then
\\[5pt]
\begin{enumerate*}[label = (\roman*), itemjoin={\\[5pt]\noindent}, mode=unboxed]
\item
The degree of each vertex $K$ in the graph of the triangulation $\graph$ is $2$.

\item
A radial triangulation has at least three elements.

\item
If $\# \triang = j$ then $\graph = C_{j}$ \emph{i.e.} it is a cycle graph of $j$ vertices.

\item
The triangulation is bipartite if and only if the number of triangles is even.
\end{enumerate*}
\begin{proof}
\begin{enumerate*}[label = (\roman*), itemjoin={\\[5pt]\noindent}, mode=unboxed]
\item
Let $\pole$ be the pole of the triangulation and $K$ an element of the triangulation, since $\pole\in \partial K$ this implies the other two vertices of $K$ must be exterior, thus $\deg(K)<3$. Since $\graph$ is connected we know $\deg(K)\geq 1$. However if $\deg (K) = 1$ this would imply that only one side of the triangle is in the interior of $\Omega$ and the three vertices of $K$ would be exterior which can not be since $\triang$ is radial. Thus, $\deg (K) = 2$ and the proof of the first part.

\item
If a triangulation has less than three elements all the vertices are exterior, therefore it can not be radial.

\item
Since $\deg(K) = 2$ for all $K\in \triang$ then $\# \triang  = \# \edges$. On the other hand, by construction $\graph$ is simple, therefore it must be the cycle graph $C_{j}$.

\item
Since $\graph$ is a cycle graph it contains a unique cycle of length $\# \triang$ then due to theorem \ref{Th bipartite graphs characterization} it is bipartite if and only if $\# \triang$ is even.
\end{enumerate*}
\end{proof}
\end{proposition}
\begin{definition}\label{Def radial subgraph}
Let $\mathcal{O}$ be an open bounded polygonal domain, $\{K: K\in \triang\}$ a triangulation of $\mathcal{O}$ and $\xi$ an interior vertex of the triangulation. We define its associated radial subgraph $C_{\xi}$ defined by the triangles $V_{\xi} \defining \{K\in \triang: \xi \in \partial K\}$ and the set of edges $\E_{\xi}$ which connect two elements of $V_{\xi}$. Clearly $C_{\xi}$ is radial and due to proposition \ref{Th radial mesh} we know $C_{\xi} = C_{j}$ where $j \defining \# V_{\xi} = \# \{K\in \triang: \xi \in \partial K\} $.
\end{definition}
\begin{lemma}\label{Th cycle and exterior vertices relationship}
Let $\{K: K\in \triang\}$ be a triangulation of $\Omega$ then, for any cycle $C$ in the graph $\graph$ and for each $K\in C$ there exists at least one vertex $\xi\in \vertexint\cap \partial K$ such that the winding number $ n(\gamma_{\scriptscriptstyle C}, \xi) $ is not zero.
\begin{proof} Since $C$ is a cycle the path $\gamma_{\scriptscriptstyle C}$ divides the plane in two components. Let $K$ be an element of the cycle, since $K$ is a triangle at least one of its vertices lies within the bounded component of the plane, namely $\xi$. Then $ n(\gamma_{\scriptscriptstyle C}, \xi) \neq 0$ and $\xi$ can not be an exterior vertex because $\Omega$ is simply connected and $n(\gamma_{\scriptscriptstyle C}, \b) = 0$ for all $\b\in \R^{\! 2} - \Omega$, see \cite{ConwayComplex}.
\end{proof}
\end{lemma}
\begin{theorem}\label{Th first inductive step}
Let $\{K: K\in \triang\}$ be a triangulation of $\Omega$ such that it has only one interior vertex, namely $\vertexint = \{\pole\}$. Then the triangulation is bipartite if and only if $\# \{K\in \triang: \pole \in \partial K\} $ is even.
\begin{proof} First notice that the triangulation has at least three elements and, as discussed in remark \ref{Rem Graph Connectedness} its graph $\graph$ is connected. If $ \{K\in \triang: \pole\in \partial K\} = \triang$ there is nothing to prove due to proposition \ref{Th radial mesh}.

If $\triang - \{K\in \triang: \pole\in \partial K\}\neq \emptyset$ consider the radial subgraph $C_{\pole}$ given by definition \ref{Def radial subgraph}, thus $C_{\pole} = C_{j}$ where $j\defining  \# \{K\in \triang: \pole\in \partial K\}$. Now let $L$ be in $\triang - \{K\in \triang: \pole\in \partial K\}$ then, it must hold that its three vertices are exterior and due to lemma \ref{Th cycle and exterior vertices relationship} $L$ can not belong to any cycle. Therefore, any cycle $C$ in the graph $\graph$ must lie in the subgraph $C_{\xi_{p}}$, however the subgraph $C_{\xi_{p}} = C_{j}$ is cyclic and its unique cycle is itself. Thus, the graph $\graph$ contains a unique cycle and due to theorem \ref{Th bipartite graphs characterization} it is bipartite if and only if $j = \# \{K\in \triang: \pole \in \partial K\} $ is even.
\end{proof}
\end{theorem}
Finally, we characterize the bipartite grids in the following result.
\begin{theorem}\label{Th bipartite grids characterization}
A triangulation $\{K: K\in \triang\}$ of $\Omega$ is bipartite if and only if for every interior vertex $\xi\in \vertexint$ the number of incident triangles $\# \{K\in \triang: \xi \in \partial K\} $ is even.
\begin{proof} We begin proving the necessity. Suppose there exists an interior vertex, namely\ $\pole\in \vertexint$ such that $\# \{K\in \triang: \pole \in \partial K\}  = 2j+1$. Then the cycle $C = \langle K_{1}, \ldots , K_{2j+1}, K_{1}\rangle$ where $\pole\in K_{i}$ for all $i = 1, \ldots, 2 j + 1$ belongs to the graph and it has odd length. Therefore the graph can not be bipartite due to theorem \ref{Th bipartite graphs characterization}.

In order to prove the sufficiency of the condition we proceed by induction on the number of interior vertices. The case $\# \vertexint  = 0$ implies that the graph $\graph$ has no cycles due to lemma \ref{Th cycle and exterior vertices relationship} and therefore it is bipartite according to theorem \ref{Th bipartite graphs characterization}. If $\# \vertexint = 1$ the result follows due to theorem \ref{Th first inductive step}. Assume now the result holds whenever the number of interior vertices is less or equal than $j$ and let $\{K:K\in \triang\}$ be a triangulation of $\Omega$ such that $\# \vertexint  = j+1$. Define
\begin{subequations}
\begin{equation*}
\triang_{\scriptscriptstyle SC} \defining \{K\in \triang: \#(\partial K\cap \vertexext) < 3\}
\end{equation*}
%
and
\begin{equation*}
\Omega_{\scriptscriptstyle SC} \defining \bigcup\{K: K\in\triang_{\scriptscriptstyle SC}\}
\end{equation*}
\end{subequations}
i.e. the elements of the triangulation which do not have three exterior vertices and
the natural subdomain of $\Omega$ for which $\triang_{\scriptscriptstyle SC}$ is a triangulation. Clearly the number of interior vertices in $\triang_{\scriptscriptstyle SC}$ equals the number of interior vertices in $\triang$ i.e. $j + 1$. Moreover, if $\xi$ is an interior vertex of the triangulation then
\begin{equation*}
\# \{K\in \triang_{\scriptscriptstyle SC}: \xi \in \partial K\}  =
\# \{K\in \triang: \xi \in \partial K\}
\end{equation*}
Recall that if $\#(\partial K\cap \vertexext) = 3$ then, due to lemma \ref{Th cycle and exterior vertices relationship}, $K$ can not belong to any cycle. Therefore, every cycle in $\triang_{\scriptscriptstyle SC}$ is even if and only if every cycle in $\triang$ is even. Hence, without loss of generality it can be assumed that the triangulation satisfies $\#(\partial K\cap \vertexext) < 3$ for all $K\in \triang$.

Let $K$ be an exterior triangle, then, one of its vertices, namely $\zeta$ must be interior and the other two exterior since only one of its sides lies on the boundary of $\Omega$ i.e. it has degree two. From now on we denote $K_{\zeta}$ this element of $\triang$.
Consider the domain $\Omega_{\zeta} \defining \Omega - cl (K_{\zeta})$, clearly the domain must be simply connected since $K_{\zeta}$ is simply connected and an exterior element; also $\zeta\notin \Omega_{\zeta}$. On the other hand, the family $\triang_{\zeta}\defining\{K\in \triang: K\neq K_{\zeta}\}$ is clearly a triangulation of the domain $\Omega_{\zeta}$ in which $\zeta$ is not an interior vertex. Hence,  $\triang_{\zeta}$ is a triangulation of a simply connected domain whose interior vertices are given by the set $\vertexint - \{\zeta\}$ i.e. it has only $j$ interior vertices; additionally  we have
\begin{equation*}
\# \{K\in \triang_{\zeta}: \xi \in \partial K\}  =
\# \{K\in \triang: \xi \in \partial K\}
\,,\quad \forall  \, \xi\in \vertexint - \{\zeta\} .
\end{equation*}
We conclude that each interior vertex of the triangulation $\triang_{\zeta}$ has an even number of incident triangles, due to the induction hypothesis the graph is bipartite; denote  $U_{\zeta}, W_{\zeta}$ the vertex bipartition of the graph $\mathcal{G}_{\scriptscriptstyle \triang_{\zeta}}$.

Consider now the radial subgraph $C_{\zeta}$ given by definition \ref{Def radial subgraph}. From the hypothesis we know $C_{\zeta}$ contains an even number of triangles, then the set $\{K\in \triang_{\zeta}: \zeta\in \partial K\}$ has an odd number of triangles. Since $\deg(K_{\zeta}) = 2$
denote $L, M\in \triang$ the triangles such that $\vert \partial K_{\zeta}\cap \partial L\vert >0$,  $\vert \partial K_{\zeta}\cap \partial M\vert >0$ clearly $K, L\in \triang_{\zeta}$, we claim these triangles belong to only one subset of the vertex bipartition.  Let $\langle L, \ldots, M\rangle$ be the unique path from $L$ to $M$ within both graphs $C_{\zeta}$ and $\mathcal{G}_{\scriptscriptstyle \triang_{\zeta}}$ evidently its length is the even number $
\# C_{\zeta}  - 2$. Recalling the graph $\mathcal{G}_{\scriptscriptstyle \triang_{\zeta}}$ is bipartite we conclude that $L$ and $M$ must belong to the same element of the vertex partition, either $U_{\zeta}$ or $W_{\zeta}$ without loss of generality assume $L, M\in U_{\zeta}$. Thus, the pair $U\defining U_{\zeta}$, $W \defining W_{\zeta}\cup \{K_{\zeta}\}$ is a vertex bipartition of the graph $\graph$ since it only has one extra element: $K_{\zeta}$ whose only two edges have the other endpoint on a triangle belonging to $U$. This completes the proof.
\end{proof}
\end{theorem}
%
%
%
%
\subsection{The Refinement}\label{Sec refinement in 2 d}
%
%
In this section we present a result of existence for bipartite triangulations of simply connected polygonal regions in $\R^{\!2}$; it is a method for refining a given grid into a bipartite one. It does not pursue advantages from the numerical point of view, only from the analytical point of view in order to make possible the mixed-mixed variational formulation of the porous media problem.
%
%
\begin{figure}[!]
\caption[1]{Bipartite Refinement}\label{Fig Grid Refinements}
\centerline{\resizebox{10cm}{8cm}
{\includegraphics{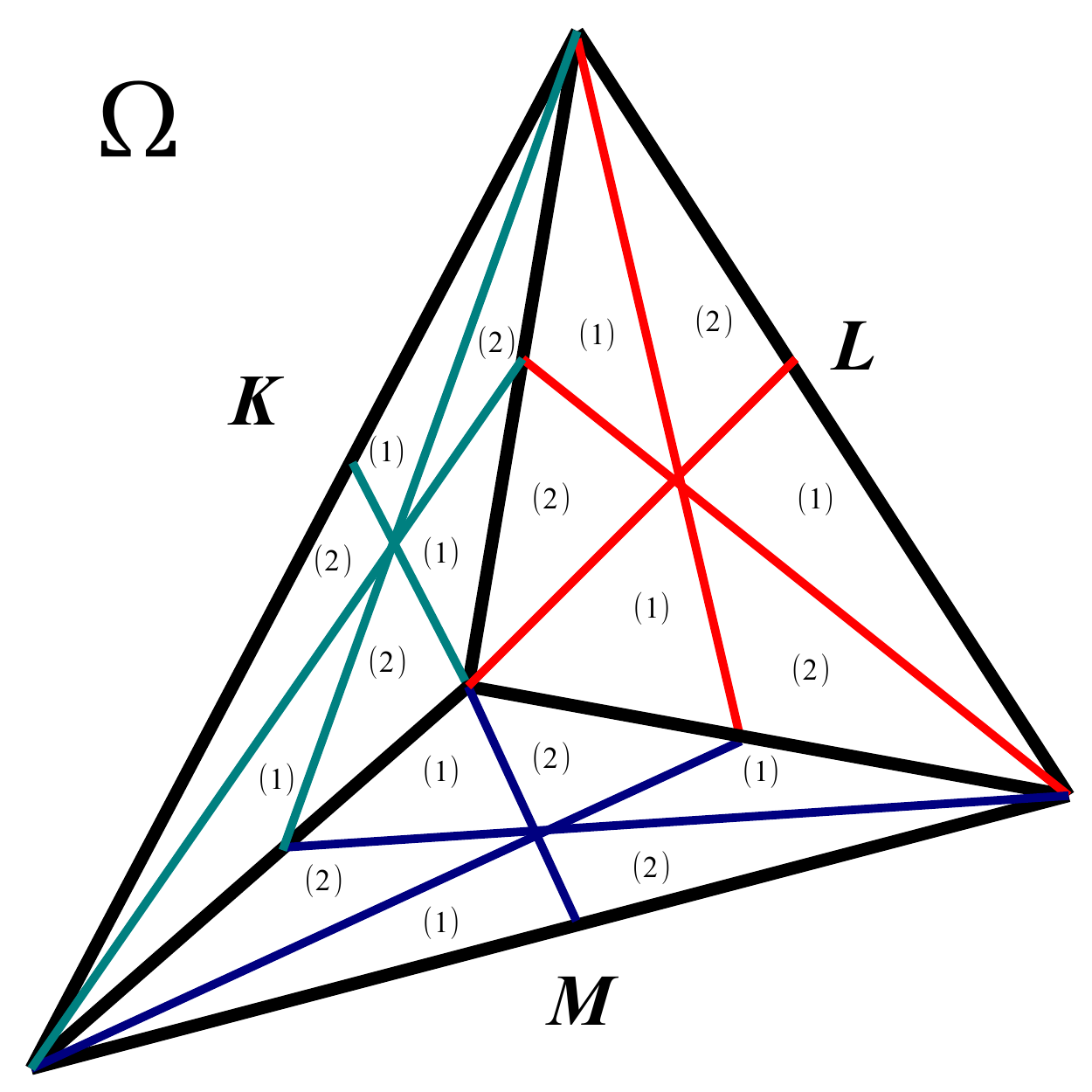}}}
\end{figure}
%
%
%
\begin{definition}\label{Def Bipartite Gridding}
\begin{enumerate*}[label = (\roman*), itemjoin={\\[5pt]\noindent}, mode=unboxed]
\item
$K\subseteq \R^{\!2}$ be a triangle, we define its BIPARTITE GRIDDING as the collection of the six triangles $\{L_{\,i}: 1\leq i\leq 6\}$ generated by the three medians of the triangle $K$.

\item
Let $\triang = \{K: K\in \triang\}$ be a triangulation of the domain $\Omega$.
We define its BIPARTITE REFINEMENT as the mesh that is generated applying the bipartite gridding process to each triangle $K$ of the triangulation. We denote it $\B\triang  = \{L: L\in \B\triang\}$.

\end{enumerate*}
\end{definition}
%
%
%
%
\begin{theorem}\label{Th gravity refinement bipartite property}
Let $\Omega\subseteq \R^{\!2}$ then it has bipartite triangulation.
%
%
%
%
\begin{proof} Let $\triang = \{K:K\in \triang\}$ be any triangulation of $\Omega$, since the domain is polygonal such triangulation exists; let $\B \triang = \{L: L\in \B\triang\}$ be its bipartite refinement; we are to prove that this grid is bipartite. Let $\xi$ be an interior vertex of the triangulation $\mathcal{B}\triang$, due to theorem \ref{Th bipartite grids characterization} it is enough to show that $\# \{L \in \B\triang: \xi\in \partial L\}$ is even. Notice that $\xi$ has only three possibilities:
\\[5pt]
\begin{enumerate*}[label = (\roman*), itemjoin={\\[5pt]\noindent}, mode=unboxed]
\item
If $\xi$ was an interior vertex of the triangulation $\triang$ then $\# \{K\in \triang: \xi\in \partial K\} \neq 0$. For each triangle $K\in \triang$ incident in $\xi$ the median passing through $\xi$ divides $K$ in two sub-triangles of $\B\triang$, thus $\# \{L\in \B\triang: \xi\in \partial L\} = 2 \# \{K\in \triang: \xi\in \partial K\} $.

\item
If $\xi = \b_{K}$ i.e it is the barycenter of some $K\in \triang$ then $\# \{L\in \B\triang: \xi\in \partial L\} = 6$.

\item
If $\xi$ is the middle point of a triangle's edge $K$ in the triangulation $\triang$, recalling it is an interior vertex we conclude that $\# \{L\in \B\triang: \xi\in \partial L\} = 4$.
\end{enumerate*}
\\[5pt]
Since in the three cases the number of triangles incident in $\xi$ is even the result follows.
\end{proof}
\end{theorem}
\begin{remark}\label{Rem refining grids}
Figure \ref{Fig Grid Refinements} illustrates the theorem above, clearly the triangulation $\triang \defining \{K, L, M\}$ of the domain $\Omega$ is not bipartite however $\B\triang$ is bipartite mesh.
\end{remark}
%
%
%
%
\begin{figure}[!]
\caption[1]{Red Refinement}\label{Fig Red Refinement}
\centerline{\resizebox{10cm}{8cm}
{\includegraphics{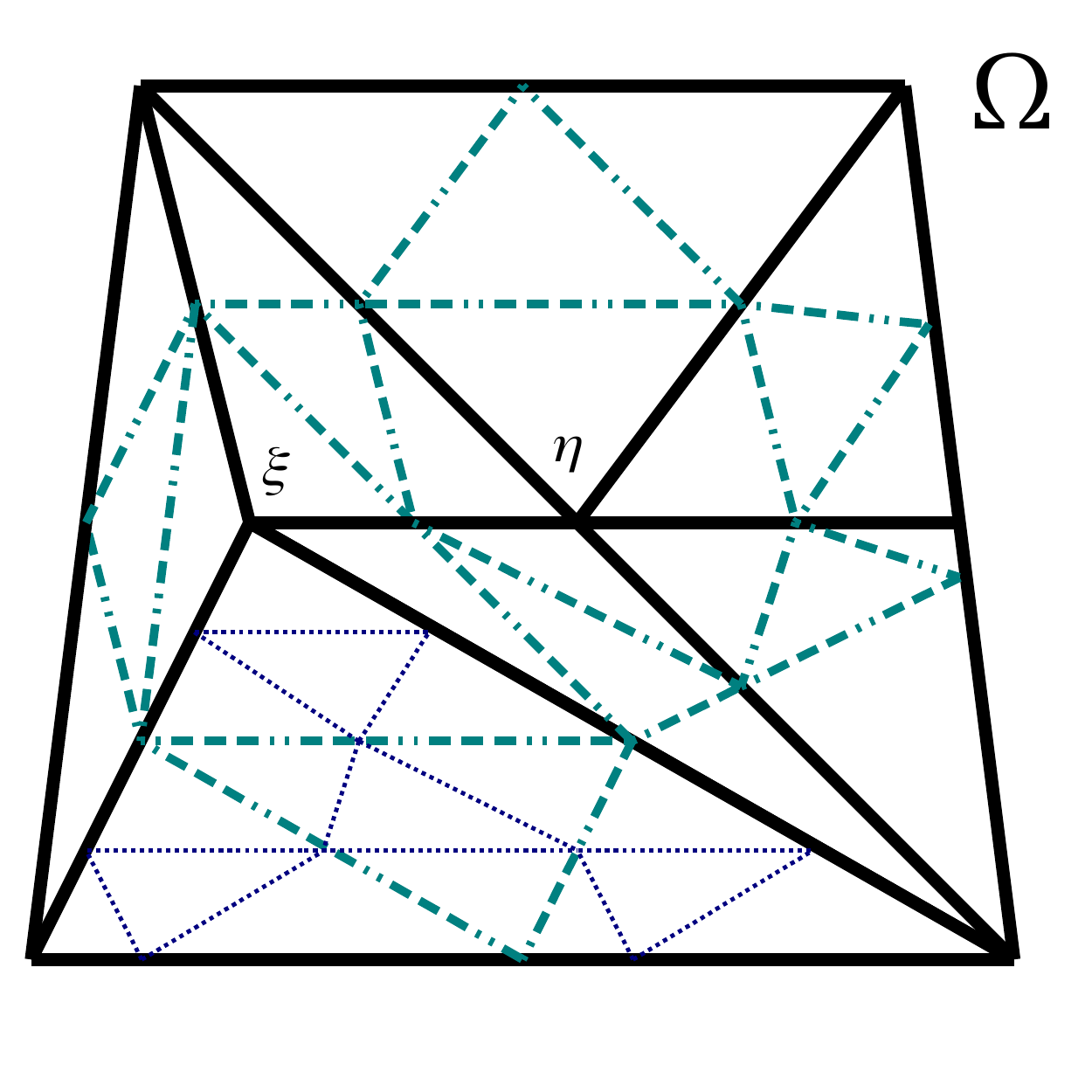}}}
\end{figure}
%
%
%
%
\subsection{Bipartite Grids and Red Refinement  }\label{Sec self-similar refinement}
%
%
In this section we propose a method to refine bipartite grids without deteriorating the regularity of the mesh and preserving the bipartite property. We start recalling the definition of the classical shape and size parameters of a grid \cite{GiraultRaviartNS}.
\begin{definition}\label{Def parameters of triangle regularity}
Let $\triang$ be a triangulation of a domain $\Omega$ and $K \in \triang$ one of its triangular elements
\\[5pt]
\begin{enumerate*}[label = (\roman*), itemjoin={\\[5pt]\noindent}, mode=unboxed]
\item
Denote $h_{\scriptscriptstyle K}$ its diameter.

\item
Denote $\rho_{\scriptscriptstyle K} \defining \sup\{ \text{diameter of } B: B \; \text{ball contained in}\; K\}$.

\item
Define the regularity of the triangle by the ratio $\zeta_{\scriptscriptstyle K}\defining \dfrac{h_{\scriptscriptstyle K}}{\rho_{\scriptscriptstyle K}}$.

\item
The size of the mesh $h$ is defined by $h\defining \max\{h_{\scriptscriptstyle K}: K\in \triang\}$.

\item
The regularity of the mesh $\triang$ is defined by $\zeta\defining \zeta(\triang) = \max\{\dfrac{h_{\scriptscriptstyle K}}{\rho_{\scriptscriptstyle K}}: K\in \triang\}$.
\end{enumerate*}
\end{definition}
Next we recall the definition of red refinement in two dimensions, see \cite{BankShermanWeiser}
\begin{definition}\label{Def Sel-similar Refinement}
\begin{enumerate*}[label = (\roman*), itemjoin={\\[5pt]\noindent}, mode=unboxed]
\item
Let $K$ be a triangle. We subdivide it into four geometrically similar triangles by pairwise connecting the midpoints of the three edges $K$. We call this family the red refinement of $K$ and denote it $\mathcal{S}K$.

\item
Let $\triang$ be a triangulation of a domain $\Omega$ we say its red refinement is the mesh obtained applying the red refinement to each of its elements and denote it $\mathcal{S} \triang$.

\item
Define recursively $\mathcal S^{(j+1)}\triang$ as the red refinement of the grid $\mathcal{S}^{(j)}\triang$. It is understood that $ \mathcal{S}^{(1)}\triang \defining \mathcal {S}\triang$.
\end{enumerate*}
\end{definition}
\begin{remark}\label{Rem self-similar refinement and regularity}
\begin{enumerate*}[label = (\roman*), itemjoin={\\[5pt]\noindent}, mode=unboxed]
\item
Notice that the size of the self-similar refinement is half the size of the original triangulation, \emph{i.e.} $h(\mathcal{S}\triang) = \dfrac{1}{2}\,h(\triang)$, see \cite{BraessFEM}.

\item
Since the red refinement of a mesh $\mathcal{S}(\triang)$ divides each triangle in four subtriangles, each of them similar to the original one, the regularity of the mesh remains equal i.e. $\zeta(\triang) = \zeta(\mathcal{S}\triang) = \zeta(\mathcal{S}^{(j)}\triang)$ for all $j\in \N$, see \cite{BraessFEM}. This is seen in figure \ref{Fig Red Refinement} which depicts two levels of red refinement for a given grid. The first one in dashed-doted line and a second level, performed only on a portion, in doted line.

\end{enumerate*}
\end{remark}
\begin{theorem}\label{Th bipartite self-similar refinement}
Let $\triang = \{K\in \triang\}$ be a bipartite triangulation of the domain $\Omega$ then its red refinement $\mathcal{S}(\triang)$ is a bipartite grid.
\begin{proof} Let $\xi$ be an interior vertex of the triangulation $\mathcal{S}\triang$, due to theorem \ref{Th bipartite grids characterization} we need to show that $\# \{L \in \mathcal{S}\triang: \xi\in \partial L\}$ is even. Notice that $\xi$ has only two possibilities:
\\[5pt]
\begin{enumerate*}[label = (\roman*), itemjoin={\\[5pt]\noindent}, mode=unboxed]
\item
$\xi$ is an interior vertex of the original triangulation. Since the self-similar refinement does not introduce new edges incident on the vertices of the original triangulation $\triang$ then
\begin{equation*}
\#\{L\in \triang:\xi\in\partial L\} = \#\{K\in \triang:\xi\in\partial K\}
\end{equation*}
Since $\triang$ is bipartite the cardinal of the right hand side is even.

\item
There exists $\sigma\in \edges$, the set of edges of the original triangulation $\triang$ such that $\xi$ is the middle point of $\sigma$. Let $K_{1}, K_{2}\in \triang$ the only pair such that $\sigma = \partial K_{1}\cap \partial K_{2}$. Then, exactly six triangles of $\mathcal{S}\triang$ are incident on $\xi$: three contained in $K_{1}$ and three contained in $K_{2}$; i.e. $\#\{L\in \mathcal{S}\triang: \xi\in \partial L \}$ is even.
\end{enumerate*}
\end{proof}
\end{theorem}
\begin{remark}\label{Rem Neutrality of self-similar refinement}
It is important to observe that the red refinement does not generate bipartite grids. In figure \ref{Fig Red Refinement} two internal vertices of the original triangulation $\xi$ and $\eta$ are depicted. Notice that, in both cases the red refinement does not change the number of triangles incident on them. Therefore, it is necessary to use the refinement given in definition \ref{Def Refining Grid} to generate bipartite grids; this is unfortunate due to the quality deterioration such procedure introduces in the mesh. However, in the practical case when bipartite grids of arbitrary small size are necessary, due to theorem \ref{Th bipartite self-similar refinement} the bipartite refinement needs to be applied only once \emph{i.e.} the grid $\mathcal{S}^{(j)} \B\triang$ for $j\in \N$ large enough satisfies b requirements.
\end{remark}
%
%
%
%
\section{The Three Dimensional Case}
%
%
We start this section recalling basic standard results for simply connected spaces \cite{Munkres}, and defining the main tools to analyze the bipartite tetrahedral grids.
\begin{definition}\label{Def topological paths and homotopy}
Let $X$ be a topological space then
\\[5pt]
\begin{enumerate*}[label = (\roman*), itemjoin={\\[5pt]\noindent}, mode=unboxed]
\item
Given two points $x, y$ of the space $X$, a \textbf{path} in $X$ from $x$ to $y$ is a continuous map $\gamma: [0, 1]\rightarrow X$ such that $\gamma(0) = x$ and $\gamma(1) = y$.

\item
The space is said to be \textbf{path connected} if every pair of points of $X$ can be joined by a path in $X$.

\item
Two paths $\gamma, \gamma\prime$ mapping the interval $I = [0,1]$ into $X$ are said to be \textbf{path homotopic} if they have the same initial point $x_{0}$ and the same final point $x_{1}$, and if there is a continuous map $F: I\times I \rightarrow X$ such that
\begin{subequations}\label{Eq homotopy definition}
\begin{align}
& F(s, 0) = \gamma(s) & & \text{and} & & F(s,1) = \gamma\prime(s)\,, & &\forall\,s\in I\\
& F(0, t) = x_{0} & & \text{and} & & F(1,t) = x_{1}\,, & &\forall\,t\in I
\end{align}
\end{subequations}
We say $\gamma$ and $\gamma\prime$ are \textbf{homotopy related} or simply \textbf{homotopic} and denote the relation by $\gamma \rel \gamma\prime$. Finally, we call $F$ a \textbf{path homotopy} between $\gamma$ and $\gamma\prime$.
\end{enumerate*}
\end{definition}
\begin{lemma}\label{Th equivalence relation}
The relation $\rel$ is an equivalence relation.
\begin{proof} See \cite{Munkres}.
\end{proof}
\end{lemma}
\begin{definition}\label{Def product of paths}
Let $X$ be a topological space, if $\gamma$ is a path in $X$ from $x_{0}$ to $x_{1}$ and if $\gamma\prime$ is a path in $X$ from $x_{1}$ to $x_{2}$ we define the \textbf{product} $\gamma*\gamma\prime$ of $\gamma$ and $\gamma\prime$ to be the path given by
\begin{equation*}
\gamma*\gamma\prime\,(s) \defining
\begin{cases}
\gamma (2s) & s\in[0, \frac{1}{2}],\\
\gamma\prime(2s-1) & s\in [0, \frac{1}{2}].
\end{cases}
\end{equation*}
\end{definition}
\begin{definition}\label{Def fundamental group and loops}
Let $X$ be a topological space then
\\[5pt]
\begin{enumerate*}[label = (\roman*), itemjoin={\\[5pt]\noindent}, mode=unboxed]
\item
Given a point $x_{0}\in X$ a path in $X$ that begins and ends at $x_{0}$ is called a \textbf{loop} based at $x_{0}$.

\item
The set of path homotopy classes of loops based at $x_{0}$ with the operation $*$, is called the \textbf{fundamental group} of $X$ relative to the \textbf{base point} $x_{0}$. It is denoted by $\pi_{1}(X, x_{0})$.

\item
The space is said to be \textbf{simply connected} if it is a path-connected space and $\pi_{1}(X, x_{0})$ is the trivial (one-element) group for some $x_{0}\in X$, and hence for every $x_{0}\in X$. We express the fact that $\pi_{1}(X, x_{0})$ is the trivial group by writing $\pi_{1}(X, x_{0}) = 0$.
\end{enumerate*}
\end{definition}
Next we introduce the basic topological spaces for the current problem
\begin{definition}\label{Def basic topological space}
Let $\mathcal{S}^{1}\subseteq \R^{\!2} = \{(x, y)\in \R^{\!2}: \vert (x, y)\vert = 1\}$ and $\jmath(\mathcal{S}^{1}) = \{(x, y, z)\in \R^{3}: \vert (x, y)\vert = 1, z = 0\}$ i.e. its ``natural" embedding in $\R^{\!3}$. Define
\begin{equation}\label{Eq basic topological space}
\bspace \defining  \R^{\!3} - \jmath(\mathcal{S}^{1}) 
%
\end{equation}
\end{definition}
\begin{theorem}\label{Th basic space simply conneccted}
The space $\bspace$ is not simply connected.
\begin{proof}
See proposition 6.1 \cite{Massey}.
\end{proof}
\end{theorem}
\begin{proposition}\label{Th homomorphism to basic space}
Let $K\subseteq \R^{3}$ be a tetrahedron, $\Phi$ one of its faces and $\mathscr{S}$ the contour of $\Phi$ \emph{i.e.} $\mathscr{S} = \partial (\Phi)$ in the trace topology of $\Phi$. Then
\\[5pt]
\begin{enumerate*}[label = (\roman*), itemjoin={\\[5pt]\noindent}, mode=unboxed]
\item
$\R^{3} - \mathscr{S}$ is homeomorphic to $\R^{3} - \jmath(\mathcal{S}^{1}) = \bspace$.

\item
$\R^{3} - \mathscr{S}$ is not simply connected.

\end{enumerate*}
\begin{proof}
\begin{enumerate*}[label = (\roman*), itemjoin={\\[5pt]\noindent}, mode=unboxed]
\item
Let $h: \R^{\!3}\rightarrow \R^{\!3}$ be any homeomorphism such that $h(\mathscr{S}) = \jmath(\mathcal{S}^{1})$ then $h\vert_{ \R^{3} - \mathscr{S}}$ is a homeomorphism.

\item
Since $\R^{\!3} - \jmath(\mathcal{S}^{1}) = \bspace$ is not simply connected as shown in theorem \ref{Th basic space simply conneccted} and its homeomorphic to $\R^{\!3} - \mathscr{S}$ the result follows.
\end{enumerate*}
\end{proof}
\end{proposition}
%
%
\subsection{Characterization of Bipartite Tetrahedral Grids}
\begin{figure}[!]
\caption[1]{Tent--type Tetrahedral Mesh}\label{Fig Tent2}
\centerline{\resizebox{10cm}{8cm}
{\includegraphics{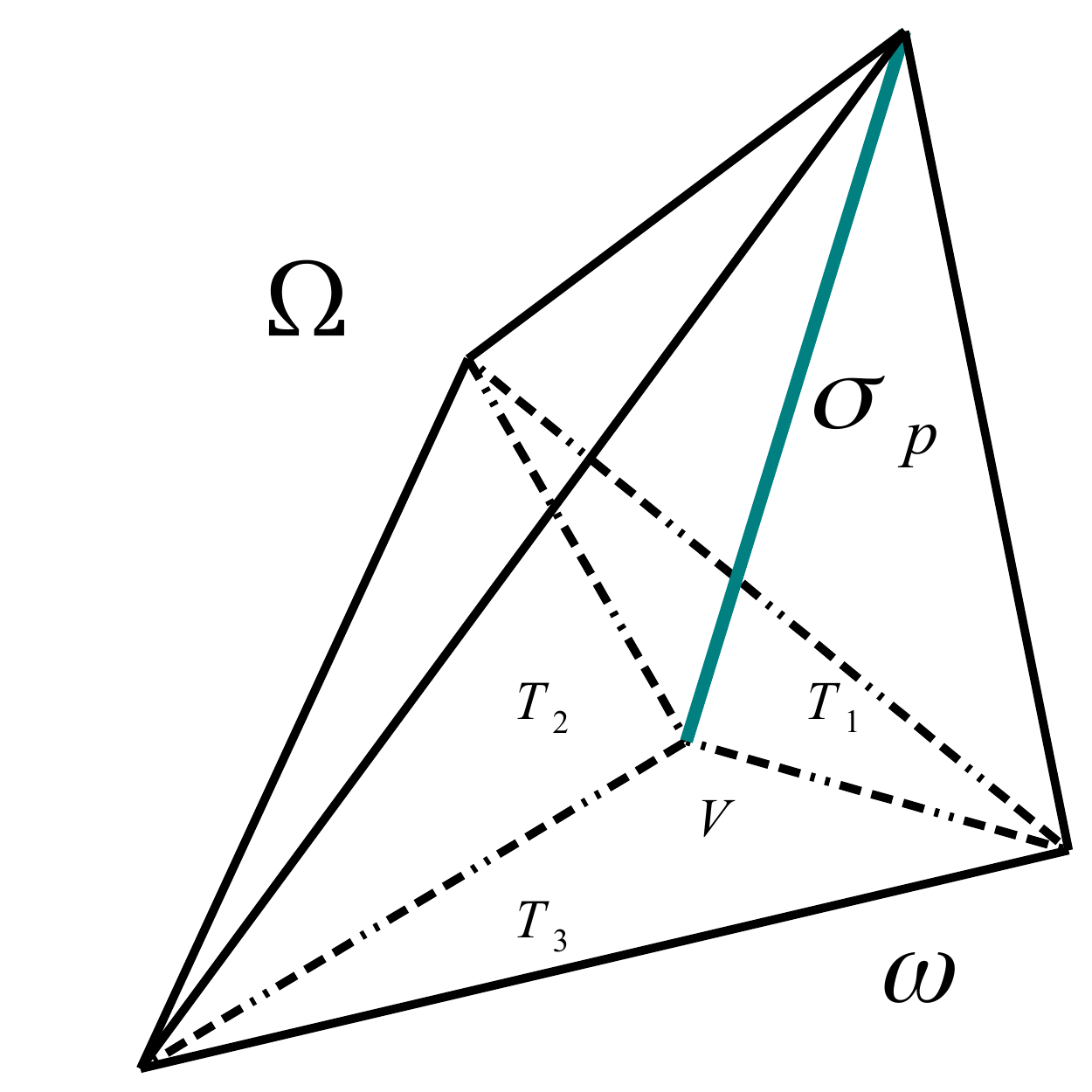}}}
\end{figure}

%
%
We start giving the definitions in order to model a tetrahedral mesh with a graph.
\begin{definition}\label{Def triangulation 3 d}
Let $\mathcal{O}\subseteq \R^{\!3}$ be an open bounded polyhedral domain and $\{K: K\in \triang\}$ be a tetrahedral mesh of $\mathcal{O}$. We denote $\vertex$, $\edges$ and $\faces$ the set of vertices, edges and faces of the tetrahedral mesh respectively.
\\[5pt]
\begin{enumerate*}[label = (\roman*), itemjoin={\\[5pt]\noindent}, mode=unboxed]

\item
We say a face $\Phi \in \faces$ is EXTERIOR if it lies on the boundary of the domain $\partial \mathcal{O}$ (if $\vert \Phi\cap  \partial \mathcal{O}\vert_{2}>0$) and denote the set of exterior faces by $\facesext$. We say a face $\Phi\in \faces$ of a tetrahedron is INTERIOR if it belongs to the interior of the domain (if $\vert \Phi\cap  \partial \mathcal{O}\vert_{2}=0$) and denote the set of interior faces by $\facesint$.

\item
The associated graph $\graph(\triang, \facesint)$ is defined in the following way. The set of vertices $\triang$ is defined by the set of tetrahedra and there is an arc in $\facesint$ between two elements $K, L\in \triang$ if $\vert \partial K\cap \partial L\vert_{2}  > 0$ i.e. if they share a common face.

\item
An element of the mesh $\{K:K\in \triang\}$ is said to be isolated if it shares no common face with any other tetrahedron. Equivalently if its degree as vertex of $\graph$ is zero. We denote $\triangcon$ the connected elements of the mesh.

\item
We say a tetrahedron is INTERIOR if $\vert \partial K\cap \partial \mathcal{O}\vert_{2} = 0$ (equivalently, if non of its faces is on the boundary of the domain) and denote $\triangint$ the set of interior tetrahedra. We say a tetrahedron $K\in \triang$ is EXTERIOR if $\vert \partial K\cap \partial \mathcal{O}\vert_{2}>0$ (equivalently, if one or more of its faces is on the boundary of the domain) and denote $\triangext$ the set of exterior tetrahedra.

\item
We say an edge $\sigma\in \edges$ is EXTERIOR if it lies on the boundary of the domain $\partial \mathcal{O}$ and denote the set of exterior edges by $\edgesext$. We say an edge $\sigma\in \edges$ of a tetrahedron is INTERIOR if it belongs to the interior of the domain and denote the set of interior edges by $\edgesint$.

\end{enumerate*}
\end{definition}
\begin{definition}\label{Def Bipartite triangulation 3 d}
Let $\mathcal{O}\subseteq \R^{\! 3}$ be an open bounded polyhedral domain and $\{K: K\in \triang\}$ be a tetrahedral mesh of $\mathcal{O}$. We say a tetrahedral mesh $\{K : K\in \triang\}$ is BIPARTITE if it can be colored with only two colors. Equivalently if its associated graph $\graph(\triang, \facesint)$ is bipartite.
\end{definition}
\begin{definition}\label{Def control points and paths 3 d}
Let $\mathcal{O}\subseteq \R^{\!3}$ be an open bounded polygonal domain and $\{K: K\in \triang\}$ be a tetrahedral mesh of $\mathcal{O}$.
\\[5pt]
\begin{enumerate*}[label = (\roman*), itemjoin={\\[5pt]\noindent}, mode=unboxed]
\item
For each $K\in \triang$ we denote $\b_{K}$ its barycenter or center of gravity.

\item
Given a cycle $C = \langle K_{1}, K_{2}, \ldots, K_{j}\rangle$ in the graph $\graph$ define its associated loop $\gamma_{\scriptscriptstyle C}$ by the sequence of segments $[\b_{\scriptscriptstyle K_{1}}, \b_{\scriptscriptstyle K_{2}}], [\b_{\scriptscriptstyle K_{2}}, \b_{\scriptscriptstyle K_{3}}]\ldots [\b_{\scriptscriptstyle K_{j-1}},\b_{\scriptscriptstyle K_{j}}]$ i.e. $\gamma_{\scriptscriptstyle C}$ is contained in $\mathcal{O}$.

\end{enumerate*}
\end{definition}
We are aimed to give necessary and sufficient conditions for a tetrahedral mesh of a domain $\Omega$ to be bipartite. First we focus on a very particular type of tetrahedral mesh.
\begin{definition}\label{Def radial triangulation 3 d}
We say a tetrahedral mesh $\{K: K\in \triang\}$ of a polyhedral domain $\mathcal{O}$ is a TENT if it is connected and has only one interior edge named POLE, denoted $\Pole$ such that $\{K\in \triang: \Pole\subseteq \partial K\} = \triang$. We call TENT GRAPH its associated graph $\graph = \mathcal{G}_{\sigma}$.
\end{definition}
\begin{remark}\label{Rem tent}
\begin{enumerate*}[label = (\roman*), itemjoin={\\[5pt]\noindent}, mode=unboxed]
\item
The figure \ref{Fig Tent2} displays the most basic case of a tent for a polyhedral (tetrahedral in this case) domain $\Omega$. A polygonal domain $\omega$ is subdivided in three triangles $T_{1}, T_{2}, T_{3}$ with only one interior vertex $V$. The pole of the mesh $\Pole$ stands on the vertex $V$ ``lifting" each triangle into a tetrahedron.

\item
A more general type of tent consists on a polygonal domain triangulated with a radial triangulation (seen in definition \ref{Def radial triangulation}) i.e. all the triangles are incident in one single vertex $V$ and the pole $\Pole$ stands on the vertex ``lifting" each triangle into a tetrahedron.

\item
In the most general version of a tent there may not exist a plane hosting one face of each tetrahedron of the mesh.
\end{enumerate*}
\end{remark}
%
%
%
%
\begin{proposition}\label{Th radial mesh 3 d}
Let $\{K: K\in \triang\}$ be a tent mesh of the domain $\Omega$ then
\\[5pt]
\begin{enumerate*}[label = (\roman*), itemjoin={\\[5pt]\noindent}, mode=unboxed]
\item
The degree of each vertex $K$ in the graph of the mesh $\graph$ is $2$.

\item
A tent mesh has at least three elements.

\item
If $\# \triang = j$ then $\graph = C_{j}$ \emph{i.e.} it is a cycle graph.

\item
The tent mesh is bipartite if and only if the number of tetrahedra is even.
\end{enumerate*}
\begin{proof}
\begin{enumerate*}[label = (\roman*), itemjoin={\\[5pt]\noindent}, mode=unboxed]
\item
Let $\Pole$ be the pole of the tent and $K$ an element of the mesh, since $\Pole\in \partial K$ and only two faces of $K$ host $\Pole$ this implies the other two faces of $K$ must be exterior, thus $\deg(K)\leq 2$. Since $\graph$ is connected we know $\deg(K)\geq 1$. However if $\deg (K) = 1$ this would imply that only one face of the tetrahedron is in the interior of $\Omega$ and the three faces of $K$ would be exterior which can not be since $\Pole$ is an interior edge; this proves the first part.

\item
If a mesh has less than three tetrahedra, i.e. two or one, there exists at most one interior face. Therefore, no interior edges exist and, by definition, the mesh can not be a tent.

\item
The fact $\deg(K) = 2$ for all $K\in \triang$ implies $\# \triang  = \# \facesint$. On the other hand, by construction $\graph$ is simple, therefore it must be the cycle graph $C_{j}$.

\item
Since $\graph$ is a cycle graph it contains a unique cycle of length $\# \triang$ then due to theorem \ref{Th bipartite graphs characterization} it is bipartite if and only if $\# \triang$ is even.
\end{enumerate*}
\end{proof}
\end{proposition}
\begin{definition}\label{Def radial subgraph 3 d}
Let $\mathcal{O}\subseteq \R^{\!3}$ be an open bounded polyhedral domain, $\{K: K\in \triang\}$ a tetrahedral mesh of $\mathcal{O}$ and $\sigma$ an interior edge. We define its associated subgraph $C_{\sigma}$ by the set of tetrahedra $V_{\sigma} \defining \{K\in \triang: \sigma \in \partial K\}$ and the set of faces $\F_{\sigma}$ which connect two elements of $V_{\sigma}$. Clearly $C_{\sigma}$ is a tent and due to proposition \ref{Th radial mesh} we know $C_{\sigma} = C_{j}$ where $j \defining \# V_{\sigma} = \# \{K\in \triang: \sigma \in \partial K\} $.
\end{definition}
\begin{lemma}\label{Th cycle and exterior vertices relationship 3 d}
Let $\{K: K\in \triang\}$ be a tetrahedral mesh of $\Omega$ then, for any cycle $C$ in the graph $\graph$ and for each $K\in C$ there exists at least one interior edge $\sigma\in \edgesint$ such that $\sigma\subseteq \partial K$.
\begin{proof} Let $C$ be a cycle and $\gamma_{\scriptscriptstyle C}$ its associated loop. Fix $K\in C$, and let $L_{1}, L_{2}\in C$ the elements such that the sequence $\langle  L_{1}, K,  L_{2} \rangle$ is in the cycle. Let $\Phi \defining \partial L_{1}\cap \partial K$ i.e. the unique face shared by $L_{1}$ and $K$, since $\vert \gamma_{\scriptscriptstyle C}\cap (L_{1}\cap K) \vert_{1} >0$
the loop hits both sides of the face $\Phi$. Let
\begin{equation*}
\mathscr {S}\defining \text{cl} \bigcup_{\sigma\,\in\, \edges} \{\sigma: \sigma\in \partial K\cap \partial L_{1}\}
\end{equation*}
i.e. the $\mathscr {S}$ is the contour of the face $\Phi$. Thus $\mathscr {S}$ is a one-dimensional manifold in $\R^{3}$ with the shape of a triangle. If $K$ has no interior edges it would hold that $\mathscr {S} \cap \Omega = \emptyset$ this would imply that $\Omega$ and consequently $\gamma_{\scriptscriptstyle C}$ are contained in $\R^{\!3} - \mathscr {S}$. However, the loop $\gamma_{\scriptscriptstyle C}$ hits both sides of the surface $\Phi$ and $\R^{\!3} - \mathscr {S}$ is not simply connected as seen in proposition \ref{Th homomorphism to basic space}, therefore $\gamma_{\scriptscriptstyle C}$ can not be homotopic to one point. This contradicts the hypothesis for $\Omega$ been simply connected. Therefore, one of the edges $\sigma$ of $K$ must be interior.
\end{proof}
\end{lemma}
\begin{theorem}\label{Th first inductive step 3 d}
Let $\{K: K\in \triang\}$ be a tetrahedral mesh of $\Omega$ such that it has only one interior edge, namely $\edgesint = \{\Pole\}$. Then the mesh is bipartite if and only if $\# \{K\in \triang: \Pole \subseteq \partial K\} $ is even.
\begin{proof} First observe that the mesh has at least three tetrahedra, other wise no edge could be interior. If $ \{K\in \triang: \Pole\in \partial K\} = \triang$ there is nothing to prove due to proposition \ref{Th radial mesh 3 d}.

If $\triang - \{K\in \triang: \pole\in \partial K\}\neq \emptyset$ consider the tent subgraph $C_{\Pole}$ given by definition \ref{Def radial subgraph 3 d}, which is cyclic i.e. $C_{\Pole} = C_{j}$ where $j\defining  \# \{K\in \triang: \Pole\in \partial K\}$.

Let $L$ be in $\triang - \{K\in \triang: \Pole\in \partial K\}$ then, it must hold that its six edges are exterior and due to lemma \ref{Th cycle and exterior vertices relationship} $L$ can not belong to any cycle. Hence, any cycle $C$ in the graph $\graph$ must lie in the subgraph $C_{\Pole}$, however the subgraph $\mathcal{G}_{\Pole} = C_{j}$ is cyclic and its unique cycle is itself. Therefore, the graph $\graph$ contains a unique cycle and due to theorem \ref{Th bipartite graphs characterization} is bipartite if and only if $j = \# \{K\in \triang: \Pole \in \partial K\} $ is even.
\end{proof}
\end{theorem}
Before characterizing the bipartite grids we need an intermediate lemma regarding the process of ``removing" an exterior tetrahedron from a triangulation.
\begin{lemma}\label{Th removing tetrahedron lemma}
Let $\{L: L\in \triang\}$ be a tetrahedral mesh of $\Omega$, and $K$ be an exterior tetrahedron of $\triang$ such that has at least one interior edge $\sigma$. Define $\triang_{K} \defining \{L\in \triang: L\neq K\}$ and the domain $\Omega_{K}\defining\Omega - \text{cl } (K)$.
\\[5pt]
\begin{enumerate*}[label = (\roman*), itemjoin={\\[5pt]\noindent}, mode=unboxed]
\item
Then the degree $\deg (K)\in\{2,3\}$ and $K$ has at most three interior edges.

\item
If $\deg(K) = 2$ then the graph $\mathcal{G}_{\scriptscriptstyle \triang_{K}}$ is connected and consequently the domain $\Omega_{K}$ is also connected.

\item
If $\deg(K) = 3$ and $K$ has two or three interior edges then the graph $\mathcal{G}_{\scriptscriptstyle \triang_{K}}$ is connected and consequently the domain $\Omega_{K}$  is also connected.

\item
If $\deg(K) = 2$ or $\deg(K) = 3$ and it has two or three interior edges the domain $\Omega_{K}$ is simply connected.

\item\label{Th case of one interior edge}
If $\deg(K) = 2$ or $\deg(K) = 3$ and has only one interior edge the domain $\Omega_{K}$, can have one or two components. However each component of $\Omega_{K}$ is simply connected.
\end{enumerate*}
\begin{proof}
\begin{enumerate*}[label = (\roman*), itemjoin={\\[5pt]\noindent}, mode=unboxed]
\item
$K$ has at least one interior edge therefore its degree has to be greater than one, additionally $K$ is exterior, then its degree has to be less than 4 i.e. $\deg(K)\in \{2,3\}$. Seeing that $K$ is exterior, it has a face $\Phi$ contained on the boundary of $\Omega$, therefore the three edges belonging to $\partial \Phi \cap \partial K$ have to be exterior, i.e. the element has at most three interior edges.

\item
Let $L_{1}, L_{2}\in \triang$ the only two tetrahedra such that $\vert \partial L_{i}\cap \partial K\vert_{2} > 0$, $i = 1, 2$. It also holds that $\sigma\subseteq \partial L_{1}\cap \partial L_{2}$ and by definition $L_{1}, L_{2}\in \triang_{ K}$. Since $\sigma\in \edgesint$ consider the tent subgraph $C_{\sigma}$ which is contained in $\graph$. Clearly $L_{1}, L_{2}\in C_{\sigma}$ and there exists a unique path $\langle L_{1}, \ldots, L_{2}\rangle$ which is contained in the cycle graph $C_{\sigma}$ but does not hit $K$. Observing that all the elements of the path belong to $\mathcal{G}_{\scriptscriptstyle \triang_{K}}$, we conclude that $L_{1}$ and $L_{2}$ can be connected in $\mathcal{G}_{\scriptscriptstyle \triang_{K}}$, therefore this graph is connected. The connectedness of $\Omega_{K}$ follows immediately.

\item
Let $L_{1}, L_{2}, L_{3}\in \triang$ the only three tetrahedra such that $\vert \partial L_{i}\cap \partial K\vert_{2} > 0$, $i = 1, 2, 3$. And let $\sigma, \tau$ be two interior edges of $K$, without loss of generality assume
\begin{equation*}
\sigma\subseteq \partial L_{1}\cap \partial L_{2}\,,\quad
\tau\subseteq \partial L_{2}\cap \partial L_{3}
\end{equation*}
By definition $L_{1}, L_{2}, L_{3}\in \triang_{ K}$. Using $C_{\sigma}$, $C_{\tau}$ and repeating the previous argument we conclude that $L_{1}$, $L_{2}$ can be connected and $L_{2}$, $L_{3}$ can be connected, consequently $L_{1}$ and $L_{3}$ can also be connected. Then, the graph $\mathcal{G}_{\scriptscriptstyle \triang_{K}}$ is connected. From here connectedness of $\Omega_{K}$ follows.

\item
The domain $\Omega_{K}$ is connected as seen in the (ii) and (iii), recalling that $\Omega$ is simply connected and that $K$ is convex and exterior we conclude $\Omega - cl (K)$ must also be simply connected.

\item
If $\deg(K) = 2$ let $L_{1}, L_{2}$ the tetrahedra connected to $K$. Then, it should hold that $\sigma \subseteq \partial L_{1} \cap \partial L_{2}$. From the second part we know $\mathcal{G}_{\scriptscriptstyle \triang_{K}}$ is connected. Additionally we know that $\Omega$ is simply connected and $K$ is convex and exterior, then we conclude $\Omega - cl(K)$ is simply connected.

If $\deg(K) = 3$ let $L_{1}, L_{2}, L_{3}$ the tetrahedra neighboring $K$ then without loss of generality we can assume $\sigma\subseteq \partial L_{1}\cap \partial L_{2}$ and due to the previous analysis $L_{1}$, $L_{2}$ can be connected. Hence $\Omega_{K}$ has at most two components. If $\Omega_{K}$ is connected it must be simply connected since $K$ is convex exterior and $\Omega$ is simply connected. If $\Omega_{K}$ is not connected then $L_{1}, L_{2}$ belong to one component, namely $\Theta_{1}$ which is simply connected by the previous argument and $L_{3}$ belongs to the other component, namely $\Theta_{2}$. Thus, $K$ is convex and exterior to the polyhedral domain $int[\Theta_{2}\cup cl(K)]\subseteq \Omega$. Recalling $\Omega$ is simply connected, we conclude that $\Theta_{2}$ is simply connected.
\end{enumerate*}
\end{proof}
\end{lemma}

Finally, we characterize the bipartite grids in the following result.
\begin{theorem}\label{Th bipartite grids characterization 3 d}
Let $\{K: K\in \triang\}$ be a tetrahedral mesh of $\Omega$, the grid is bipartite if and only if for every interior edge $\sigma\in \edgesint$ the number of incident tetrahedra $\# \{K\in \triang: \sigma \subseteq \partial K\} $ is even.
\begin{proof} We begin proving the necessity. Suppose there exists an interior edge, namely\ $\Pole\in \edgesint$ such that $\# \{K\in \triang: \Pole \subseteq\partial K\}  = 2j+1$; considering the tent subgraph $\mathcal{G}_{\Pole}$ and its unique cycle $C = \langle K_{1}, K_{2}, \ldots , K_{2j+1}, K_{1}\rangle$ we conclude that the mesh can not be bipartite due to theorem \ref{Th bipartite graphs characterization}.

In order to prove the sufficiency of the condition we proceed by induction on the number of interior edges. The case $\# \edgesint  = 0$ implies that the graph $\graph$ has no cycles due to lemma \ref{Th cycle and exterior vertices relationship 3 d} and therefore it is bipartite according to theorem \ref{Th bipartite graphs characterization}. If $\# \edgesint = 1$ the result follows due to theorem \ref{Th first inductive step 3 d}. Assume now the result holds whenever the number of interior edges is less or equal than $j$. Define
\begin{subequations}
\begin{equation*}
\triang_{\scriptscriptstyle SC} \defining \{K\in \triang: \#(\partial K\cap \edgesext) < 6\}
\end{equation*}
and
\begin{equation*}
\Omega_{\scriptscriptstyle SC} \defining \bigcup\{K: K\in\triang_{\scriptscriptstyle SC}\} .
\end{equation*}
\end{subequations}
i.e. the tetrahedra of the mesh which do not have six exterior edges, and the natural subdomain of $\Omega$ for which $\triang_{\scriptscriptstyle SC}$ is a triangulation. Clearly the number of interior edges in $\triang_{\scriptscriptstyle SC}$ equals the number of interior edges in $\triang$ i.e. $j + 1$. Moreover, if $\sigma$ is an interior edge of the tetrahedral mesh then
\begin{equation*}
\# \{K\in \triang_{\scriptscriptstyle SC}: \sigma \subseteq \partial K\}  =
\# \{K\in \triang: \sigma \subseteq  \partial K\}
\end{equation*}
Recall that if $\#(\partial K\cap \edgesext) = 6$ then due to lemma \ref{Th cycle and exterior vertices relationship 3 d} $K$ can not belong to any cycle, consequently the cycles in $\triang_{\scriptscriptstyle SC}$ and in $\triang$ are the same. Moreover, a cycle is even in $\triang_{\scriptscriptstyle SC}$ if and only if is even in $\triang$. Hence, without loss of generality it can be assumed that the tetrahedral grid satisfies $\#(\partial K\cap \edgesext) < 6$ for all $K\in \triang$.

Let $K$ be an exterior tetrahedron, then one of its edges, namely $\sigma$ must be interior. From now on we denote $K_{\sigma}$ this element of $\triang$. Consider the domain $\Omega_{\,\sigma} \defining \Omega - cl (K_{\sigma})$, by definition $\sigma \nsubseteq \Omega_{\sigma}$. On the other hand the family $\triang_{\sigma}\defining\{K\in \triang: K\neq K_{\sigma}\}$ is clearly a tetrahedral mesh of the domain $\Omega_{\,\sigma}$ where $\sigma$ is not an interior edge. Since $K_{\sigma}$ is exterior but has an interior edge then, due to lemma \ref{Th removing tetrahedron lemma} its degree is two or three and it could have one, two or three interior edges one of which must be $\sigma$. In any of the cases the set of interior edges $\edgesint(\triang_{\sigma})$ of the triangulation $\triang_{\sigma}$ is contained in $\edgesint - \{\sigma\}$ i.e. it has at most $j$ interior edges. Moreover, the fact that only one exterior tetrahedron was removed from $\triang$ yields
\begin{equation*}
\# \{K\in \triang_{\sigma}: \tau \subseteq \partial K\}  =
\# \{K\in \triang: \tau \subseteq  \partial K\}
\,,\quad \forall  \, \tau \in \edgesint(\triang_{\sigma})
\end{equation*}
Therefore, due to the hypothesis, each interior edge of the tetrahedral mesh $\triang_{\sigma}$ has an even number of incident tetrahedra. Before the induction hypothesis can be applied several cases have to be analyzed.
\\[5pt]
\begin{enumerate*}[label = \bf \Alph*., itemjoin={\\[5pt]\noindent}, mode=unboxed]
\item
\underline{$\deg (K_{\sigma}) = 2$} or
\underline{$\deg(K_{\sigma}) = 3$ with $\#(\edgesint \cap K_{\sigma})\in\{2, 3\}$}. In this case, due to lemma \ref{Th removing tetrahedron lemma} the domain $\Omega_{\sigma}$ is simply connected and since it has at most $j$ interior edges the induction hypothesis implies the graph $\mathcal{G}_{\scriptscriptstyle \triang_{\sigma}}$ is bipartite, i.e. there exists a vertex bipartition of the graph $U_{\sigma}, W_{\sigma}$. Next we analyze all the possible subcases.
\begin{asparaenum}[(i)]
\item\label{Case degree two and one interior edge}
\noindent \underline{$\deg(K_{\sigma}) = 2$ and $\#\{\tau\in \edgesint: \sigma\subseteq \partial K\} = \#\{\sigma\} = 1$}. Consider the tent subgraph $C_{\sigma}$ given by definition \ref{Def radial subgraph}. From the hypothesis we know $C_{\sigma}$ contains an even number of tetrahedra, then the set $\{K\in \triang_{\sigma}: \sigma\subseteq \partial K\}$ has an odd number of tetrahedra. Let $L_{1}, L_{2}\in \triang$ be the two tetrahedra such that $\vert \partial K_{\sigma}\cap \partial L\vert_{i} >0$,  $\sigma\subseteq \partial L_{i}$ for $i = 1,2$. By definition $L_{1}, L_{2}\in \triang_{\sigma}$, we claim these tetrahedra belong to only one subset of the vertex bipartition.  Let $\langle L_{1}, \ldots, L_{2}\rangle$ be the unique path from $L_{1}$ to $L_{2}$ within both graphs $C_{\sigma}$ and $\mathcal{G}_{\scriptscriptstyle \triang_{\sigma}}$ clearly, it has even length $\# C_{\sigma}  - 2$. Then since the graph $\mathcal{G}_{\scriptscriptstyle \triang_{\sigma}}$ is bipartite $L_{1}$ and $L_{2}$ belong to the same element of the vertex partition, either $U_{\sigma}$ or $W_{\sigma}$ without loss of generality assume $L_{1}, L_{2}\in U_{\sigma}$. Thus, the pair $U\defining U_{\sigma}$, $W \defining W_{\sigma}\cup \{K_{\sigma}\}$ is a vertex bipartition of the graph $\graph$ since it only has one extra element $K_{\sigma}$, whose only two edges have the other endpoint on a tetrahedron belonging to $U$. The proof is complete for this case.

\item\label{Case two interior edges}
\noindent \underline{$\#\{\tau\in \edgesint: \sigma\subseteq \partial K\} = \#\{\sigma, \tau\} = 2$}. In this case it must hold  $\deg(K_{\sigma}) = 3$. Let $L_{1}, L_{2}, L_{3}$ the tetrahedra such that $\vert\partial L_{i}\cap \partial K\vert_{2} > 0$ for $i = 1, 2, 3$, without loss of generality we can assume
\begin{equation}\label{Eq incidence on edges}
\sigma\subseteq \partial L_{1}\cap \partial L_{2}\,,\quad \tau\subseteq \partial L_{2}\cap \partial L_{3}
\end{equation}
Notice that $L_{1}$, $L_{2}$ and $L_{3}$ belong to the graphs $\mathcal{G}_{\scriptscriptstyle \triang_{\sigma}}$. Let $C_{\sigma}$ and $C_{\tau}$ the tent subgraphs given by definition \ref{Def radial subgraph} then both have an even number of tetrahedra. Also, due to \eqref{Eq incidence on edges} $L_{1}, L_{2}\in C_{\sigma}$ and $L_{2}, L_{3}\in C_{\tau}$ must hold. Let $\langle L_{1}, \ldots, L_{2}\rangle$ be the unique path from $L_{1}$ to $L_{2}$ within both graphs $C_{\sigma}$ and $\mathcal{G}_{\scriptscriptstyle \triang_{\sigma}}$ clearly, it has even length $\# C_{\sigma}  - 2$. Then, since the graph $\mathcal{G}_{\scriptscriptstyle \triang_{\sigma}}$ is bipartite, $L_{1}$ and $L_{2}$ must belong to the same subset of the vertex partition, either $U_{\sigma}$ or $W_{\sigma}$. Repeating the argument $L_{2}$ and $L_{3}$ must belong to the same subset either $U_{\sigma}$ or $W_{\sigma}$, therefore we conclude the three of them belong to one single set, without loss of generality assume $L_{1}, L_{2}, L_{3}\in U_{\sigma}$. Hence, the pair $U\defining U_{\sigma}$, $W \defining W_{\sigma}\cup \{K_{\sigma}\}$ is a vertex bipartition of the graph $\graph$ since it only has one extra element $K_{\sigma}$, whose only three edges have the other endpoint on a tetrahedron belonging to $U$. The case has been proved.

\item \label{Case three interior edges}
\noindent \underline{$\#\{\tau\in \edgesint: \sigma\subseteq \partial K\} = \#\{\sigma, \tau, \varrho\} = 3$}. In this case it must hold  $\deg(K_{\sigma}) = 3$. Let $L_{1}, L_{2}, L_{3}$ the tetrahedra such that $\vert\partial L_{i}\cap \partial K\vert_{2} > 0$ for $i = 1, 2, 3$ and such that
\begin{equation}\label{Eq incidence on edges case 2}
\sigma\subseteq \partial L_{1}\cap \partial L_{2}\,,\quad \tau\subseteq \partial L_{2}\cap \partial L_{3}
\,,\quad \varrho\subseteq \partial L_{3}\cap \partial L_{1}
\end{equation}
This case is reducible to the previous one \ref{Case two interior edges} since it is enough to analyze the tent subgraphs $C_{\sigma}$ and $C_{\tau}$ to conclude $L_{1}, L_{2}$ and $L_{3}$ belong to the same subset of the bipartition $U_{\sigma}$ of the graph $\mathcal{G}_{\scriptscriptstyle \triang_{\sigma}}$. Using the same arguments as before we have that the pair $U\defining U_{\sigma}$, $W \defining W_{\sigma}\cup \{K_{\sigma}\}$ is a vertex bipartition of the graph $\graph$.
\end{asparaenum}

\item
\underline{$\deg(K_{\sigma}) = 3$ and $\#(\edgesint \cap \partial K_{\sigma}) = 1$}. Let $L_{1}, L_{2}, L_{3}$ the three neighboring tetrahedra to $K$ and assume that $\sigma\subseteq \partial L_{1}\cap \partial L_{2}$. In this case, due to lemma \ref{Th removing tetrahedron lemma}\ref{Th case of one interior edge} the domain $\Omega_{\sigma}$ has one or two connected two components. If it has only one component the problem is reduced to the case \ref{Case three interior edges}.

If $\Omega_{\sigma}$ has two components, namely $\Theta_{1}, \Theta_{2}$ due to lemma \ref{Th removing tetrahedron lemma}\ref{Th case of one interior edge} each of them is simply connected and $L_{1}, L_{2}\in \Theta_{1}$, $L_{3}\in \Theta_{2}$. Using the induction hypothesis, both graphs $\mathcal{G}_{\scriptscriptstyle \Theta_{1}}$ and $\mathcal{G}_{\scriptscriptstyle \Theta_{2}}$ are bipartite. Let $(U_{\sigma}, W_{\sigma})$ and $(U', W')$ be the vertex bipartition pairs of both graphs where $L_{1}, L_{2}\in U_{\sigma}$ and $L_{3}\in U'$. We will prove that the pair $U\defining U_{\sigma}\cup U'$, $W\defining W_{\sigma}\cup \{K_{\sigma}\}\cup W'$ is a vertex bipartition for the graph $\graph$. Observe that if a closed path $C$ hits elements on both subgraphs $\mathcal{G}_{\Theta_{1}}, \mathcal{G}_{\Theta_{2}}$ it can not be a cycle. Let $C$ be a closed path, since the one element that connects the subgraphs $\mathcal{G}_{\Theta_{1}}, \mathcal{G}_{\Theta_{2}}$ is $K_{\sigma}$, then it must be part of the cycle $C$. However, there is no way of ``crossing" from $\mathcal{G}_{\Theta_{1}}$ to $ \mathcal{G}_{\Theta_{2}}$ or viceversa in a closed path without having to hit $K_{\sigma}$ twice, therefore $C$ can not be a cycle. Now let $C$ be a cycle in the graph $\graph$, due to the previous discussion only two cases are possible.
\begin{asparaenum}[(i)]
\item
\noindent \underline{$C$ belongs to the subgraph of the domain $\Theta_{1}\cup K_{\sigma}$}.
This case is reducible to the case \ref{Case degree two and one interior edge} since $K_{\sigma}$ has degree two and only one interior edge in the subgraph corresponding to the mesh inherited to the domain $\Theta_{1}\cup K_{\sigma}$. Therefore, the pair $U \defining U_{\sigma}$ and $V \defining W_{\sigma} \cup \{K_{\sigma}\}$ constitutes a vertex bipartition for the graph of the subdomain $\Theta_{1} \cup \{K_{\sigma}\}$, then the length of $C$ must be even.

\item
\noindent \underline{$C$ belongs to the subgraph of the domain $\Theta_{2}\cup K_{\sigma}$}. In this case since $K_{\sigma}$ has no edges which are interior to the domain $\Theta_{2}\cup K_{\sigma}$. Then, as seen in lemma \ref{Th cycle and exterior vertices relationship 3 d} the element $K_{\sigma}$ can not belong to the cycle. Thus $C$ must be contained in the graph of the domain $\Theta_{2}$ which is bipartite, i.e. its length must be even.

\end{asparaenum}
Since in both cases the length of the cycle is even, the proof the case B is complete, this finishes the proof of the theorem.
\end{enumerate*}
\end{proof}
\end{theorem}
%
%
\subsection{Existence of Bipartite Grids }   
%
%
We finish the section presenting a result for the existence of a bipartite tetrahedral grids for simply connected polyhedral domains in $\R^{\!3}$. As in the two dimensional setting presented in section \ref{Sec refinement in 2 d} the process will deteriorate the quality of the grid, such deterioration will be estimated in section \ref{Sec refinement in 3 d}. We start recalling some well-known properties of the tetrahedron
\begin{theorem}\label{Th tetrahedron centroid properties}
Let $\Delta$ be a non-degenerate tetrahedron, and $\{\a_{\ell}: 1\leq \ell\leq 4\}$ be its four vertices, with $\a_{\ell} = (x_{\ell}, y_{\ell}, z_{\ell})$ then, the center of gravity $\overline{\a} = (\overline{x}, \overline{y}, \overline{z})$ satisfies $\overline{\a} = \frac{1}{4} \sum_{\ell\,=\,1}^{4} \a_{\ell}$.
%
%
Moreover, if $\b = (b_{x}, b_{y}, b_{z})$ is the barycenter of the face of the tetrahedron defined by the first three vertices $\{\a_{\ell}: 1\leq \ell\leq 3\}$ 
%
%
Then the points $\a_{4}, \overline{\a}$ and $\b$ are collinear. 
%
%
\end{theorem}
Recall now the standard definition of the symmetric group of a set of ``$n$ letters".
\begin{definition}\label{Def permutation of n letters}
Given $n\in \N$ we define $S_{n}$ as the set of all possible permutations of the set $\{1, 2, \ldots, n\}$.
\end{definition}
Next we define a convenient tetrahedral mesh for a non-degenerate tetrahedron
\begin{definition}\label{Def bipartite tetrahedron gridding}
Let $\Delta$ be a tetrahedron and $\{\a_{\ell}: 1\leq \ell\leq 4\}$ its vertices. For each $\pi\in S_{4}$ define
\begin{equation}\label{Eq general tetrahedron in bipartite gridding}
\Delta_{\pi}\defining \mathrm{co}\left\{\frac{1}{i}\sum_{\ell\,=\,1}^{i}\a_{\pi(\ell)}: 1\leq i\leq 4\right\}
\end{equation}
where $\mathrm{co}(A)$ denotes the convex hull of the set $A$. We define the family
\begin{equation}\label{Eq characterization bipartite gridding}
\B \Delta = \left\{\Delta_{\pi}: \pi\in S_{4}\right\}
\end{equation}
as the BIPARTITE GRIDDING of $\Delta$.
\end{definition}
\begin{remark}\label{Rem comments of barycentric gridding}
Some observations are in order for the definition above.
\\[5pt]
\begin{enumerate*}[label = (\roman*), itemjoin={\\[5pt]\noindent}, mode=unboxed]
\item
Notice that for any $\pi\in S_{4}$, the vertices of $\Delta_{\pi}$ are a vertex of $\Delta$, the midpoint of an edge of $\Delta$, the barycenter of one face of $\Delta$ and finally the centroid of $\Delta$.

\item
For any $\Delta_{\pi}$, three of its vertices are convex combinations of points belonging to extreme sets of $\Delta$: a vertex, an edge and a face. On the other hand since $\Delta$ is non-degenerate, these three vertices of $\Delta_{\pi}$ can not be collinear.

\item
Due to theorem \ref{Th tetrahedron centroid properties} the fourth vertex of $\Delta_{\pi}$ is the centroid of $\Delta$ and since the original tetrahedron is non-degenerate it can not be coplanar with the aforementioned three vertices of $\Delta_{\pi}$.

\item
Notice that $\# \B\Delta = \# S_{4} = 24$.
\end{enumerate*}
\end{remark}
%
%
%
%
\begin{theorem}\label{Th bipartite tetrahedron gridding proof}
Let $\Delta\subseteq \R^{\!3}$ be a non-degenerate tetrahedron, and $\{L_{\,i}: 1\leq i\leq 24\}$ be its bipartite gridding, then the associated graph is bipartite.
\begin{proof} Let $\{L_{\,i}: 1\leq i\leq 24\}$ be the bipartite gridding of $\Delta$ we will classify the associated vertices in the following subsets in terms of the original tetrahedron $\Delta$
\\[5pt]
\begin{enumerate*}[label = (\roman*), itemjoin={\\[5pt]\noindent}, mode=unboxed]
\item
$\{\a_{\ell}: 1\leq \ell\leq 4\}$ the vertices belonging to $\Delta$.

\item
$\xi_{\scriptscriptstyle \Delta}$ center of gravity of $\Delta$.

\item
$\{\b_{\ell}: 1\leq \ell\leq 6\}$ the barycenters of each face of $\Delta$.

\item
$\{\m_{\ell}: 1\leq \ell\leq 6\}$ the midpoints of each of the edges of $\Delta$.
\end{enumerate*}
In the bipartite gridding $\{L_{\,i}: 1\leq i\leq 24\}$ it is clear that an edge $\sigma$ is interior if and only if one of its end points is the center of gravity $\xi_{\scriptscriptstyle \Delta}$; this leaves three possible subcases:
\\[5pt]
\begin{enumerate*}[label = \alph*), itemjoin={\\[5pt]\noindent}, mode=unboxed]
\item
The other end of $\sigma$ is one of the vertices $\{\a_{\ell}: 1\leq \ell\leq 6\}$ belonging to $\Delta$. Then, six tetrahedra concur to it: two from each face concurrent to the vertex $\a_{\ell}$.

\item
The other end of $\sigma$ is one of the barycenters $\{\b_{\ell}: 1\leq\ell\leq 4\}$ of the faces of $\Delta$, namely $\Phi$. In such case, six tetrahedra concur to it, all of them having one of its faces contained in $\Phi$.

\item
The other end of $\sigma$ is one of the midpoints $\{\m_{\ell}: 1\leq\ell\leq 6\}$ of the edges of $\Delta$. In this case four tetrahedra concur to it: two from each face concurrent to the edge that hosts $\m_{\ell}$.
\end{enumerate*}
Since in all the cases the number of tetrahedra concurrent to an interior edge is even the result holds.
\end{proof}
\end{theorem}
Now we introduce a new definition
\begin{definition}\label{Def Refining Grid}
Let $\triang = \{K: K\in \triang\}$ be a tetrahedral mesh of the polyhedral domain $\Omega$.
We define its BIPARTITE REFINEMENT as the mesh that is generated applying the bipartite gridding process to each tetrahedron $K$ of the mesh. We denote it $\B\triang  = \{L: L\in \B\triang\}$.
\end{definition}
%
%
%
%
\begin{theorem}\label{Th gravity refinement bipartite property 3 d}
Let $\Omega\subseteq \R^{3}$ be a simply connected polyhedral region, then, there exists a tetrahedral mesh whose associated graph is bipartite.
\begin{proof} Let $\triang = \{K:K\in \triang\}$ be any tetrahedral mesh of $\Omega$, since the domain is polyhedral such mesh exists, and denote $\B \triang = \{L: L\in \B\triang\}$ its bipartite refinement. Let $\sigma$ be an interior edge of the tetrahedral mesh $\mathcal{B}\triang$, due to theorem \ref{Th bipartite grids characterization 3 d} we need to show that $\# \{L \in \B\triang: \sigma\in \partial L\}$ is even. Notice that $\sigma$ has only three possibilities:
\\[5pt]
\begin{enumerate*}[label = (\roman*), itemjoin={\\[5pt]\noindent}, mode=unboxed]
\item
If there exists $K\in \triang$ such that $ \sigma\subseteq \mathrm{int}(K) $ then it is one of the edges of $\B K$ . As seen in theorem \ref{Th bipartite tetrahedron gridding proof} every edge interior to one single tetrahedron generated by the bipartite gridding process has an even number of tetrahedra incident on it.

\item
If $\sigma$ was part of an interior edge of the tetrahedral mesh $\triang$ then $\# \{K\in \triang: K \; \text{is incident on } \sigma\} \neq 0$. For each tetrahedron $K\in \triang$ incident in $\sigma$ the bipartite refinement $\B\triang$ generates two tetrahedra concurrent on $\sigma$. Therefore, the number of tetrahedra incident on $\sigma$ belonging to $\B\triang$, is twice the number of tetrahedra incident on $\sigma$ belonging to $\triang$ i.e. such number is even.

\item
If there exists $K\in \triang$ such that the interior of one of its faces contains $\sigma$, then there must also exist another tetrahedron $L \in \triang$ sharing the face that contains $\sigma$. Two tetrahedra concur to this edge from $\B K$ and other two from $\B L$ giving a total of four tetrahedra of $\B \triang$ incident on $\sigma$ which is an even number.
\end{enumerate*}
\newline
In the three cases above the number of tetrahedra concurrent on $\sigma$ is even, then the result follows.
\end{proof}
\end{theorem}
%
%
\section{The Refinement}\label{Sec refinement in 3 d}
%
%
In this section we propose a method to generate bipartite grids of arbitrary small size and bounded regularity, i.e. the quality of the grids does not degenerate.
%
%
\subsection{Tetrahedron Shape Parameters}
%
%
We start this section with several definitions of geometrical shape parameters for later discussion of their relationships; we have
\begin{definition}\label{Def Tetrahedron Shape Parameters}
Let $\Delta$ be a non-degenerate tetrahedron in $\R^{\!3}$, define
\begin{subequations}\label{Eq Tetrahedron Shape Parameters}
\begin{equation}\label{Def diameter tetrahedron}
\hdelta \defining \text{diameter of the tetrahedron }\; \Delta
\end{equation}
\begin{equation}\label{Def diameter inscribed sphere}
\rhodelta \defining \sup \{\text{diameter of }\; B: B\;\text{is a ball contained in}\; \Delta\}
\end{equation}
The regularity of the tetrahedron
\begin{equation}\label{Def regularity of tetrahedron}
\regdelta \defining \frac{\hdelta}{\rhodelta}
\end{equation}
The \emph{radius ratio} is defined by
\begin{equation}\label{Def aspect ratio}
\adelta \defining 3\,\frac{r_{\mathrm{in}}}{r_{\mathrm{circ}}} 
\end{equation}
Where $r_{\mathrm{in}}, r_{\mathrm{circ}}$ are respectively, the inradius and circumradius of $\Delta$, see \cite{LiuJoe2}. Finally the \emph{mean} ratio defined in \cite{LiuJoe1} is given by
\begin{equation}\label{Def mean ratio}
\mdelta \defining \frac{12(3\,\vert\Delta\vert_{3})^{2/3}}{\sum\{\vert\sigma\vert^{2}: \sigma \text{ is an edge of }\;\Delta\}}
\end{equation}
Here $\vert \Delta\vert_{3}$ is the volume of the tetrahedron.
\end{subequations}
\end{definition}
Next we recall some previous results for relationship between shape parameters
\begin{theorem}\label{Th relation between aspect ratio and mean ratio}
For any tetrahedron $\Delta$ holds
\begin{equation}\label{Ineq aspect ratio, mean ratio equivalence}
\mdelta^{\,3} \leq  \adelta \leq \frac{2}{\sqrt[4]{6}}\, \mdelta^{ \, 3/4}
\end{equation}
Furthermore, the lower bound is optimal and tight, and the upper bound is optimal.
\begin{proof} See \cite{LiuJoe2}.
\end{proof}
\end{theorem}
On the other hand it is direct to see
\begin{equation*}
\regdelta = \frac{\hdelta}{\rhodelta} =  \frac{\hdelta}{2\,r_{\mathrm{in}}} \leq
\frac{2 \, r_{\mathrm{circ}}}{2\,r_{\mathrm{in}}}
= \frac{3}{3\,\frac{r_{\mathrm{in}}}{ r_{\mathrm{circ}}}}
\end{equation*}
Where the inequality holds since the diameter of the tetrahedron is at most the diameter of the circumradius, i.e. $\hdelta \leq 2\,r_{\mathrm{circ}}$. Then $\regdelta\leq \frac{3}{\adelta}$ and therefore
%
%
%
%
\begin{equation}\label{Ineq regularity, mean ratio control}
\regdelta\leq \frac{3}{\mdelta^{\,3}}
\end{equation}
%
%
%
%
%
Next we study the deterioration of regularity when applying the bipartite refinement to a given tetrahedron.
%
%
\subsection{Bipartite Refinement and Regularity Deterioration}
%
%
%
%
\begin{theorem}\label{Th bipartite refinement on mean ratio}
Let $K$ be a non-degenerate tetrahedron then, for any $L\in \B K$ holds
\begin{equation}\label{Ineq bipartite refinement mean ratio bound}
\frac{1}{\eta_{\scriptscriptstyle L}} \leq \frac{36\sqrt[3]{9}}{ \eta_{\scriptscriptstyle K}}
\end{equation}
\begin{proof} Fix $L\in \B K$; denote $\{\b_{i}: 1\leq i\leq 4\}$ its vertices and $\{\a_{i}: 1\leq i\leq 4\}$ be the vertices of $K$. We start assuming that the tetrahedron $L$ is the one generated by the identity permutation i.e. $L = \Delta_{I\!d}$ and its vertices satisfy
%
%
%
\begin{equation}\label{Eq pointwise notation}
\b_{i} = \frac{1}{i}\sum_{\ell\,=\,1}^{i}\a_{\ell} \,,\quad 1\leq i\leq 4
\end{equation}
%
%
Clearly $\{\vert\a_{i} - \a_{j}\vert: 1\leq i< j\leq 4\}$ and $\{\vert\b_{i} - \b_{j}\vert: 1\leq i< j\leq 4\}$ are the lengths of the six edges of $K$ and $L$ respectively, due \eqref{Eq pointwise notation} we have
\begin{equation*}
\vert \b_{1} - \b_{2}\vert =  \left\vert \a_{1} - \frac{\a_{1} + \a_{2}}{2}\right\vert
= \frac{1}{2}  \left\vert \a_{1} - \a_{2}\right\vert
\leq \frac{1}{2} \sum_{1\,\leq \,i \,< \,j \,\leq \, 4} \vert \a_{i} - \a_{j} \vert
\end{equation*}
Exhausting the remaining five cases with the same technique we have
\begin{equation*}
\vert \b_{k} - \b_{\ell}\vert
\leq \frac{1}{\lcm\{k, \ell\}} \sum_{1\,\leq \,i \,< \,j \,\leq \, 4} \vert \a_{i} - \a_{j} \vert
\leq \frac{1}{2} \sum_{1\,\leq \,i \,< \,j \,\leq \, 4} \vert \a_{i} - \a_{j} \vert, \quad 1\leq k  < \ell\leq 4
\end{equation*}
Here $\lcm\{k, \ell\}$ denotes the lowest common multiple of $k$ and $\ell$.
Therefore
\begin{multline*}
\left(\sum_{1\,\leq \,i \,< \,j \,\leq \, 4} \vert \b_{i} - \b_{j} \vert^{2}\right)^{1/2} \leq
\sqrt{6} \, \max\{\vert \b_{i} - \b_{j}\vert: 1\leq i < j\leq 4\} \\
\leq \frac{\sqrt{6}}{2} \sum_{1\,\leq \,i \,< \,j \,\leq \, 4} \vert \a_{i} - \a_{j} \vert\leq \frac{\sqrt{6}}{2}\,\sqrt{6}\left(\sum_{1\,\leq \,i \,< \,j \,\leq \, 4} \vert \a_{i} - \a_{j} \vert^{2}\right)^{1/2}
\end{multline*}
Where the factor $\sqrt{6}$ shows due to the equivalence norms $\Vert \cdot\Vert_{\infty}-\Vert \cdot\Vert_{2}$ and $\Vert \cdot\Vert_{1}-\Vert \cdot\Vert_{2}$ norms in $\R^{\!6}$; thus
\begin{equation}\label{Ineq equivalence of edges length for barycentric}
\sum_{1\,\leq \,i \,< \,j \,\leq \, 4} \vert \b_{i} - \b_{j} \vert^{2} \leq 9
\sum_{1\,\leq \,i \,< \,j \,\leq \, 4} \vert \a_{i} - \a_{j} \vert^{2}
\end{equation}
On the other hand, since the refinement $\B K = \{L_{j}: 1\leq j\leq 24\}$ is made through the centroid of $K$ all the tetrahedra have the same volume i.e. $\vert L \vert_{3} = \dfrac{1}{24}\,\vert K\vert_{3}$. Thus, recalling the definition of mean ratio given by equality \eqref{Def mean ratio} and combining it with the inequality above \eqref{Ineq equivalence of edges length for barycentric} we get
\begin{multline*}
\eta_{\scriptscriptstyle L} = \frac{12 (3 \vert  L \vert_{3} )^{2/3}}
{\sum \{\vert \b_{i} - \b_{j}\vert^{2}: 1\leq i < j \leq 4 \}} 
=
\frac{12 (3\,\frac{1}{24} \vert K\vert_{3} )^{2/3}}{\sum \{\vert \b_{i} - \b_{j}\vert^{2}: 1\leq i < j \leq 4 \}}\\
\geq
\frac{1}{9}\,\frac{1}{24^{\scriptscriptstyle 2/3}}
\frac{12 (3 \vert K \vert_{3})^{2/3}}
{\sum \{\vert \a_{i} - \a_{j}\vert^{2}: 1\leq i < j \leq 4 \}} = \frac{1}{36\sqrt[3]{9}}\;\eta_{\scriptscriptstyle K}
\end{multline*}
which is the desired estimate \eqref{Ineq bipartite refinement mean ratio bound}.

Finally, for any other tetrahedron $M\in \B K$, let $\{\c_{\,i}: 1\leq i\leq 4\}$ be its vertices, then there exists $\pi \in S_{4}$ such that
\begin{equation*}
\c_{\,i} = \frac{1}{i} \sum_{\ell\,=\,1}^{i}\a_{\pi(\ell)}\,,\quad \forall\, 1\leq i\leq 4 .
\end{equation*}
Repeating the previous arguments we have
\begin{equation*}
\sum_{1\,\leq \,i\,<\,j\,\leq 4}\vert \c_{\,i} - \c_{\,j}\vert^{2} \leq
9 \sum_{1\,\leq \,i\,<\,j\,\leq 4}\vert \a_{\pi(i)} - \a_{\pi(j)}\vert^{2}
= 9 \sum_{1\,\leq \,i\,<\,j\,\leq 4}\vert \a_{i} - \a_{j}\vert^{2}
\end{equation*}
and the estimate \eqref{Ineq bipartite refinement mean ratio bound} follows.
\end{proof}
\end{theorem}
\begin{corollary}\label{Th estimate of bipartite deterioration}
Let $K$ be a non-degenerate tetrahedron then, for any $L\in \B K$ holds
\begin{equation}\label{Ineq regularity bipartite controlled by mean ratio}
\zeta_{\scriptscriptstyle L} \leq
\frac{2^{6} \cdot  3^{9}}{\eta_{\scriptscriptstyle K}^{3}}
\end{equation}
\begin{proof} A straightforward combination of inequalities \eqref{Ineq regularity, mean ratio control} and \eqref{Ineq regularity bipartite controlled by mean ratio}.
\end{proof}
\end{corollary}
%
%
\subsection{Regularity and the QLRS Refinement } \label{Sec QLRS refinement}
%
%
We close this section citing a result given by the QLRS (quality local refinement based on subdivision) algorithm which is described in \cite{LiuJoe3}.
\begin{theorem}\label{Th regular arbitrary small bipartite grid}
Let $\triang$ be a tetrahedral mesh of the domain $\Omega$. Let $L^{(n)}$ be a refined tetrahedron produced by QLRS where $n$ denotes the number of refinement levels. Let $K\in \triang$ be the tetrahedron such that $L^{(n)}\subseteq K$, then
\begin{equation}\label{Ineq QLRS control}
\eta( L^{(n)}) \geq \frac{\sqrt[3]{4}}{11} \, \eta (K)
\end{equation}
Where $\eta( L^{(n)})$ and $\eta( K)$ stand for the mean ratio shape parameters of $L^{(n)}$ and $K$ respectively.
\end{theorem}
\begin{remark} Notice that the QLRS algorithm can be applied as many times as needed keeping a positive lower bound for its mean ratio regardless of the number of times it is applied. Therefore, using an original tetrahedral mesh this can be refined in order to generate a tetrahedral grid of prescribed size $h > 0$ and its tetrahedra have mean ratio bigger or equal than $\frac{\sqrt[3]{4}}{11} \, \eta$; where $\eta \defining \min \{\mdelta: \Delta\in \triang\} > 0$.
\end{remark}
\begin{definition}\label{Def QLRS notation}
Let $\{K:K\in\triang\}$ be a tetrahedral mesh of the domain $\Omega$, we denote by $\refin^{n} \triang$ the $n$ levels of QLRS refinement of the mesh $\triang$.
\end{definition}
\begin{theorem}\label{Th bounded regularity refinement}
Let $\triang$ be a tetrahedral mesh of $\Omega$ then the sequence $\{\B\refin^{\,n} T: n\in\N\}$ has bounded regularity, moreover there exists a positive constant $\kappa$ such that
\begin{equation}\label{Ineq bounded regularity}
\sup_{n\,\in \,\N}\max\{\zeta_{\scriptscriptstyle L}: L\in \B\refin^{\,n} T\} \leq \frac{\kappa}{\eta^{\,3}}.
\end{equation}
Where $\eta \defining \min\{\mdelta: \Delta\in \triang\}$.
\begin{proof} Let $L$ be a tetrahedron in $\B\refin^{\,n} \triang$, then there exists a $K\in \refin^{n} \triang$ and $\Delta\in \triang$ such that $L \subseteq K\subseteq \Delta$. Then combining inequalities \eqref{Ineq regularity bipartite controlled by mean ratio} and \eqref{Ineq QLRS control} we have
\begin{equation*}
\zeta_{\scriptscriptstyle L} \leq \frac{3^{6}}{\eta_{\scriptscriptstyle K}^{3}} \leq
\frac{11\cdot 2^{6} \cdot 3^{9}}{\sqrt[3]{4}} \, \frac{1}{\mdelta^{\,3}} \leq
\frac{11\cdot  2^{6} \cdot 3^{9} }{\sqrt[3]{4}} \, \frac{1}{\min \{\mdelta^{\,3}: \Delta\in \triang\}}
\end{equation*}
Therefore, the result follows for $\kappa = \dfrac{11\cdot  2^{6} \cdot 3^{9}}{\sqrt[3]{4}}$.
\end{proof}
\end{theorem}
%
%
\section{Concluding Remarks and Discussion}
%
%
The present work yields several conclusions listed below
\\[5pt]
\begin{enumerate*}[label = (\roman*), itemjoin={\\[5pt]\noindent}, mode=unboxed]
\item
We have provided a method for generating bipartite grids of prescribed size and controlled regularity in 2-D and 3-D. This provides the theoretical setting to assure the well-posed mixed-mixed formulation of a problem such as the porous media equation presented in \cite{MoralesShow2}.

\item
The method is far from been optimal. The bipartite refinements given in definitions \ref{Def Bipartite Gridding} and \ref{Def Refining Grid} for 2-D and 3-D respectively, subdivide internal angles of the initial elements. This deteriorates severely the shape quality of the mesh, as seen in the proof of theorem \ref{Th bounded regularity refinement} where the bound is amplified by a factor of $2^{6}\cdot 3^{9}$ from the constant provided by the QLRS refinement method.

\item
One initial improvement for the bipartite refinements would discuss the placement of the new internal vertex. For instance in 2-D the barycenter of the triangle could be replaced by the incenter in order to bisect the angles of the original triangle. This deteriorates the internal angles in a more balanced way.

\item
The bipartite refinement needs to be applied the the whole mesh, regardless of a-posteriori estimation of the solution or other guideline. This constraint needs to be addressed: from the mesh generation point of view (e.g. local refinement criteria or developing a different generation technique) and from the mixed-mixed variational formulation \cite{MoralesShow2} point of view.

\item
Using the bipartite refinement generates a grid which has 6 times and 24 times the number of elements of the original grid for 2-D and 3-D respectively; therefore the number of elements increases considerably. However, the computational costs for finding the center of gravity and code implementation are low, in contrast with the demanded in calculating an optimal point as the incenter previously suggested.

\end{enumerate*}
%
%
\section{Acknowledgements}
%
%
The Authors wish to thank Universidad Nacional de Colombia, Sede Medell\'in for supporting this work.
%
%

\vspace{2cc}

\bibliographystyle{plain}

\vspace{1cc}


{\small
\noindent Fernando A Morales\\
Escuela de Matem\'aticas\\
Universidad Nacional de Colombia, Sede Medell\'in\\
Calle 59 A No 63-20 - Bloque 43, oficina 106\\
Medell\'in - Colombia \\
E-mail: famoralesj@unal.edu.co\\[3pt]

\noindent Mauricio A Osorio\\
Escuela de Matem\'aticas\\
Universidad Nacional de Colombia, Sede Medell\'in\\
Calle 59 A No 63-20 - Bloque 43, oficina 106\\
Medell\'in - Colombia \\
E-mail: maosorio@unal.edu.co

}\end{document}